\documentclass[12pt, a4paper, reqno]{amsart}
\usepackage{amscd,amssymb,amsmath, array,latexsym,amsbsy}
\usepackage[arrow,matrix,curve,cmtip,ps]{xy}
\allowdisplaybreaks[2]
\begin{document}

\def\supp{\operatorname{supp}}
\def\Ind{\operatorname{Ind}}
\def\Aut{\operatorname{Aut}}
\def\Rep{\operatorname{Rep}}
\def\Res{\operatorname{Res}}

\def\Ad{\operatorname{Ad}}
\def\range{\operatorname{range}}
\def\sp{\operatorname{span}}
\def\clsp{\overline{\operatorname{span}}}
\def\dashind{\operatorname{\!-Ind}}
\def\id{\operatorname{id}}
\def\rt{\operatorname{rt}}
\def\lt{\operatorname{lt}}

\def\rtimesHoHo{\rtimes}
\def\rtimesHeHe{\rtimes}

\def\H{\mathcal{H}}
\def\L{\mathcal{L}}
\def\K{\mathcal{K}}
\def\D{\mathcal{D}}
\def\C{\mathbb{C}}
\def\R{\mathbb{R}}

\newtheorem{thm}{Theorem}  [section]
\newtheorem{cor}[thm]{Corollary}
\newtheorem{lemma}[thm]{Lemma}
\newtheorem{prop}[thm]{Proposition}
\newtheorem{goal}[thm]{Goal}
\newtheorem{thm1}{Theorem}
\theoremstyle{definition}
\newtheorem{defn}[thm]{Definition}
\newtheorem{remark}[thm]{Remark}
\newtheorem{example}[thm]{Example}
\newtheorem{remarks}[thm]{Remarks}
\newtheorem{claim}[thm]{Claim}
\newtheorem{problem}{Problem}

\numberwithin{equation}{section}

\title[Extension problems for representations]{\boldmath{Extension problems for 
representations of crossed-product $C^*$-algebras}}

\author[an Huef]{Astrid an Huef}
\address{School of Mathematics\\
The University of New South Wales\\
NSW 2052\\
Australia}
\email{astrid@unsw.edu.au}

\author[Kaliszewski]{S. Kaliszewski}
\address{Department of Mathematics\\Arizona State University\\Tempe\\ AZ 
85287-1804\\USA} \email{kaliszewski@asu.edu}

\author[Raeburn]{Iain Raeburn}
\address{School of Mathematical and Physical Sciences\\
University of Newcastle\\
NSW 2308\\ Australia}
\email{Iain.Raeburn@newcastle.edu.au}

\author[Williams]{Dana P. Williams}
\address{Department of Mathematics\\Dartmouth College\\ Hanover, NH 03755\\USA}
\email{dana.williams@dartmouth.edu}

\thanks{This research was supported by grants from the Australian Research Council, 
the National Science Foundation, the University of New South Wales and the 
Ed Shapiro Fund at Dartmouth College.}


\date{November 9, 2006}

\maketitle

\section{Introduction}

We consider the following extension problem for covariant representations of $C^*$-dynamical systems, and the analogue of this extension problem for  the dual systems involving coactions of non-abelian groups.

\begin{problem}\label{prob1}
Let $\alpha$ be an action of a locally compact group $G$ on a $C^*$-algebra $A$, let $H$ be a closed subgroup of $G$, and let $(\pi,U)$ be a covariant representation of the system $(A,H,\alpha)$ on a Hilbert space $\H$.
For which closed subgroups $K$ of $G$ containing $H$ is there a covariant representation $(\pi,V)$ of $(A,K,\alpha)$ on $\H$ such that $U=V|_H$?
\end{problem} 

When $A=\C$, $H=N$ is normal and $K=G$, this is a classical problem which has been studied using a variety of methods (see \cite{duflo,fabec,CKS,aHR-twisted}, for example), and its solution for irreducible representations is a crucial ingredient in the Mackey machine \cite{mackey}. In \cite{newPPL}, we noted that Problem~\ref{prob1} for $C^*$-dynamical systems is an interesting test question for the theory of non-abelian duality.  In particular, we showed there is a condition  on the induced representation $\Ind_N^G (\pi\rtimesHoHo U)$ of the crossed product $A\rtimesHoHo_\alpha G$ which is equivalent to the existence of $(\pi,V)$ \cite[Corollary~4]{newPPL}. For fixed $K$, we can apply this theorem to obtain a criterion involving the induced representation $\Ind_N^K(\pi\rtimesHoHo U)$ of $A\rtimesHoHo_\alpha K$. This is not a very satisfactory solution to Problem~\ref{prob1}, though, since it requires that we consider all the induced representations $\Ind_N^K(\pi\rtimesHoHo U)$ as $K$ varies.  In our first main theorem, we describe a criterion which uses the same induced representation  $\Ind_N^G(\pi\rtimesHoHo U)$ of $A\rtimesHoHo_\alpha G$ for every subgroup $K$ (see Theorem~\ref{thm-prob1}).

Our proof of Theorem~\ref{thm-prob1}  uses ideas from non-abelian duality for crossed products of $C^*$-algebras, and hence it is natural to consider also the analogue of Problem~\ref{prob1} for crossed products by coactions. In stating Problem~\ref{prob1}, we made implicit use of our ability to restrict $\alpha$ to actions of the subgroups $H$ and $K$. Coactions of $G$ restrict to coactions of quotients of $G$, and hence the most natural dual analogue of Problem~\ref{prob1} involves a pair of closed normal subgroups of $G$. A precise statement of this dual analogue is given in Problem~\ref{prob2} in \S\ref{sec-state2}. 

Our solutions to Problems~\ref{prob1} and~\ref{prob2} follow the same general pattern. Each proof has two main ingredients: a theorem describing one aspect of the duality between induction and restriction of representations, of the sort proved in \cite{KQR, EKR, BE}, and an imprimitivity theorem which allows us to recognise induced representations. When dealing with duality for crossed products of $C^*$-algebras, we have to make choices: we can use full crossed products, in which case we need to use the maximal coactions of \cite{EKQ, KQFM}, or we can use reduced crossed products, in which case we need to use the Quigg-normal coactions of \cite{Q, BE}. We prove versions of our solution to Problem~\ref{prob2} for both maximal and Quigg-normal coactions; the theorems look similar, but pose different technical problems. In solving these problems, we have proved new results which have independent interest, including a theorem about the duality of induction and restriction for maximal coactions, and a version of Green's imprimitivity theorem for reduced crossed products.

\smallskip

We begin with a short section on  induction processes and the associated imprimitivity theorems, which we hope will clarify some issues which arise in applying concrete Morita equivalences. In \S\ref{sec-prob1}, we present our solution to Problem~\ref{prob1} when $H$ is normal (see Theorem~\ref{thm-prob1}). An essential ingredient of the proof is a construction of a  Morita equivalence   which yields an induction process and an imprimitivity theorem for crossed products by dual coactions (Proposition~\ref{prop-PPL}). In Proposition~\ref{prop-cd}, we prove the duality of induction and restriction for this induction process.

In \S\ref{sec-state2}, we recall some properties of coactions, and give a detailed statement of Problem~\ref{prob2}, which is the analogue of Problem~\ref{prob1} for crossed products by coactions. In \S\ref{sectmax}, we present our solution to Problem~\ref{prob2} for maximal coactions. Here duality involves full crossed products and Green's imprimitivity theorem suffices, but we need to establish the appropriate induction-restriction result for maximal coactions (Proposition~\ref{resindmax}). At the end of \S\ref{sectmax} we indicate why our attempts to extend our results to crossed products by homogeneous spaces require that we also consider Quigg-normal coactions. 

In \S\ref{sectnor}, we solve Problem~\ref{prob2} for Quigg-normal coactions. This time the necessary induction-restriction result was established in \cite{BE}, but we need to prove a new version of Green's imprimitivity theorem for reduced crossed products (Theorem~\ref{redimpthm}). This is in itself an interesting new application of non-abelian duality: the statement is entirely about crossed products by actions, but the proof uses crossed products by coactions in a non-trivial way. We discuss an example which shows that our proof works only for normal subgroups, and observe that this example should be remembered when trying to establish universal properties of crossed products by coactions of homogeneous spaces.

In each of Sections~\ref{sec-prob1}, \ref{sectmax} and \ref{sectnor}, we use a different generalisation of the induction process and imprimitivity theorem of Mansfield for crossed products by coactions~\cite{M}. It is natural to ask to what extent these different generalisations are compatible, and we discuss this in an appendix.

\subsection*{Conventions}
Let $K$ be a closed subgroup of a locally compact group $G$. 
We use left Haar measures, and denote by $\Delta_G$ and $\Delta_K$ the modular functions on $G$ and $K$. 
If $N$ is a closed normal subgroup of $G$ with $N\subset K$,  we choose Haar measure on $K/N$ such that 
\[
\int_K f(s)\, ds=\int_{K/N}\int_N f(sn)\, dn\, d(sN)
\]
for $f\in C_c(K)$. We denote by $\lambda^G$ and $\rho^G$ the left- and right-regular representations of $G$ on $L^2(G)$. The group $G$ acts on the left of the homogeneous space $G/K$, and this induces an action $\lt:G\to\Aut C_0(G/K)$ defined by $\lt_t(f)(sK)=f(t^{-1}sK)$. When $N$ is a closed normal subgroup of $G$, there is also a right action $\rt:G/N\to\Aut C_0(G/N)$  given  by  $\rt_{tN}(f)(sN)=f(stN)$.

All homomorphisms and representations of $C^*$-algebras are assumed to be $*$-preserving. A homomorphism $\pi$ of a $C^*$-algebra $A$ into the multiplier algebra $M(B)$ of another $C^*$-algebra is \emph{non-degenerate} if $\{\pi(a)b:a\in A,b\in B\}$ is dense in $B$. A non-degenerate homomorphism $\pi:A\to M(B)$ extends uniquely to a homomorphism $\overline{\pi}:M(A)\to M(B)$, and $\overline{\pi}$ is strictly continuous; to avoid complicating formulas, we often write $\pi$ for $\overline{\pi}$. All representations of a $C^*$-algebra $A$ on Hilbert space are assumed to be non-degenerate, and we denote by $\Rep A$ the category of non-degenerate representations of $A$ on Hilbert space.
 
If $A$ and $B$ are $C^*$-algebras, a right-Hilbert $A$--$B$ bimodule is a right Hilbert $B$-module $X$ together with a homomorphism $\phi$ of $A$ into the $C^*$-algebra $\L(X)$ of adjointable operators on $X$; in practice, we suppress $\phi$ and write $a\cdot x$ for $\phi(a)x$. As in \cite{BE}, we view  a right-Hilbert $A$--$B$ bimodule $X$ as a morphism from $A$ to $B$, and say that the diagram 
\begin{equation*}
\xymatrix{
A
\ar[rr]^X
\ar[d]_Z
&&
B
\ar[d]^Y
\\
C
\ar[rr]^W
&&
D
}
\end{equation*}
commutes if the right-Hilbert $A$--$D$ bimodules $X\otimes_B Y$ and $Z\otimes_C W$ are isomorphic. Such a commuting square induces a commuting square of maps on representations:
\begin{equation*}
\xymatrix{
\Rep A
&&
\Rep B
\ar[ll]_{X\dashind_B^A}
\\
\Rep C
\ar[u]^{Z\dashind_C^A}
&&
\Rep D.
\ar[u]_{Y\dashind_D^B}
\ar[ll]_{W\dashind_D^C}
}
\end{equation*}
When $X$ is an $A$--$B$ imprimitivity bimodule, so that the map $X\dashind$ has a natural inverse $\widetilde{X}\dashind$ (see \cite[Theorem~3.29]{tfb}), we write $\cong$ beside the arrow to emphasise that it is invertible.

All tensor products of $C^*$-algebras in this paper are spatial. We consider only the full coactions of locally compact groups on $C^*$-algebras which are defined using the full group $C^*$-algebra, and our main reference for material on coactions and their crossed products is Appendix A in \cite{BE}. We assume further that all coactions are non-degenerate in the sense that $\overline{\delta(B)(1\otimes C^*(G))}=B\otimes C^*(G)$; see \cite[\S A3]{BE} for details.

\section{Induction processes and imprimitivity theorems}\label{sec-genimp}

We begin by making some general remarks about induction processes and imprimitivity theorems. Suppose that $X$ is an $A$--$B$ imprimitivity bimodule and that $\phi:C\to M(A)$ is a non-degenerate homomorphism of another $C^*$-algebra $C$ into the multiplier algebra $M(A)$ of $A$. The left action of $A$ on $X$ induces an isomorphism of $A$ onto the algebra $\K(X)$ of compact operators on $X$, which extends to an isomorphism of $M(A)$ onto the $C^*$-algebra $\L(X)$ of adjointable operators. Thus $\phi$ gives us a left action of $C$ by adjointable operators on the Hilbert $B$-module $X$, so that we can also view $X$ as a right-Hilbert $C$--$B$ bimodule. We can now apply  \cite[Theorem~3.29]{tfb} to both right-Hilbert bimodules ${}_AX_B$ and ${}_CX_B$, and thus obtain two different induction processes on representations. We are going to try to distinguish these two processes using the following convention: we denote by $X_B^A\dashind$ the functor from $\Rep B$ to $\Rep A$ given by tensoring with the imprimitivity bimodule ${}_AX_B$, and by $X\dashind_B^C$ the functor  from $\Rep B$ to $\Rep C$ given by  tensoring with the right-Hilbert bimodule ${}_CX_B$. Thus by definition we have
\[
X\dashind_B^C\pi=(X_B^A\dashind \pi)\circ \phi,
\]
where we have silently used the non-degeneracy of $X_B^A\dashind \pi$ to extend it to $M(A)$.
We may also use  abbreviations for the algebras $A, B$ and $C$ to simplify notation.

We can apply both constructions whenever we have a Morita equivalence between a crossed product $C\rtimesHoHo_{\alpha} G$ and another $C^*$-algebra $B$, when we can take for $\phi$ the canonical embedding $i_C:C\to M(C\rtimesHoHo_\alpha G)$. The resulting induction process $X\dashind_B^C$ then comes with an imprimitivity theorem:

\begin{prop}\label{genimpthm}
Suppose $\alpha:G\to \Aut C$ is an action of a locally compact group $G$ on a $C^*$-algebra $C$ and $X$ is a $(C\rtimesHoHo_\alpha G)$--$B$ imprimitivity bimodule. Let $\pi$ be a representation of $C$ on a Hilbert space $\H_\pi$. Then there is a representation $\tau$ of $B$ such that $\pi$ is unitarily equivalent to $X\dashind_B^C\tau$ if and only if there is a unitary representation $U$ of $G$ on $\H_\pi$ such that $(\pi, U)$ is a covariant representation of $(C,G,\alpha)$.
\end{prop}

\begin{proof}
Suppose $\pi:C\to B(\H_\pi)$ is unitarily equivalent to $X\dashind_B^C\tau$ for some representation $\tau$ of $B$. The isomorphism of $M(C\rtimesHoHo_\alpha G)$ onto $\L(X_B)$ induced by the left action is continuous for the strict topology on $M(C\rtimesHoHo_\alpha G)$ and the strong-operator topology on $\L(X_B)$ \cite[Proposition~C.7]{tfb}. Since $\L(X_B)$ acts non-degenerately on $X\otimes_B \H_\tau$, there is a strongly continuous representation $W$ of $G$ on $X\otimes_B \H_\tau$ such that
\[
W_s(x\otimes_B h)=(i_G(s)\cdot x)\otimes_B h\ \mbox{ for $x\in X$ and $h\in \H_\tau$}.
\]
The covariance of $(i_C,i_G)$ implies that $(X\dashind_B^C\tau,W)$ is a covariant representation of $(C,G,\alpha)$, and moving $W$ over to $\H_\pi$ gives the required unitary representation $U$.

Conversely, if there is such a representation $U$, we take 
\[
\tau:=\widetilde{X}_{C\rtimesHoHo G}^B\dashind(\pi\rtimesHoHo U),
\]
where $\widetilde{X}$ is the dual of the imprimitivity bimodule $X=X^{C\rtimesHoHo G}_B$, as in \cite[page~49]{tfb}. Since $X\otimes_B\widetilde X$ is isomorphic to ${}_{C\rtimesHoHo G}(C\rtimesHoHo_\alpha G)_{C\rtimesHoHo G}$, 
\[
X_B^{C\rtimesHoHo G}\dashind\tau=(X\otimes_B\widetilde X)\dashind(\pi\rtimesHoHo U)
\]
is equivalent to $\pi\rtimesHoHo U$, and hence
\[
X\dashind_B^{C}\tau=(X_B^{C\rtimesHoHo G}\dashind\tau)\circ i_C
\]
is equivalent to $(\pi\rtimesHoHo U)\circ i_C=\pi$.
\end{proof}

As a further illustration of this circle of ideas, we discuss Green's imprimitivity theorem. Suppose $\alpha$ is an action of a locally compact group $G$ on a $C^*$-algebra $A$, and $H$ is a closed subgroup of $G$. We denote by $X_H^G(A)$ or $X_H^G$ the $((A\otimes C_0(G/H))\rtimesHoHo_{\alpha\otimes\lt}G)$--$(A\rtimesHoHo_\alpha H)$ imprimitivity bimodule constructed by Green, which is a completion of $C_c(G, A)$ \cite[Theorem~6]{G}; we use the formulas for the  actions and inner products on dense subspaces given in \cite[Equations~B.5]{BE}. The natural embeddings 
\[
i_A:A\to A\otimes 1\subset M(A\otimes C_0(G/H))\to M((A\otimes C_0(G/H))\rtimesHoHo_{\alpha\otimes\lt}G)
\]
and $i_G:G\to UM((A\otimes C_0(G/H))\rtimesHoHo_{\alpha\otimes\lt}G)$ form a covariant representation of $(A,G,\alpha)$, and hence give a non-degenerate homomorphism
\[
i_A\rtimesHeHe i_G:A\rtimesHoHo_\alpha G\to M((A\otimes C_0(G/H))\rtimesHoHo_{\alpha\otimes\lt}G).
\]
Composing $X_H^G\dashind$ with the homomorphism
$i_A\rtimesHeHe i_G$ gives an induction process
\[
X\dashind_H^G=X\dashind_{A\rtimes H}^{A\rtimes G}:\Rep(A\rtimesHoHo_\alpha H)\to \Rep(A\rtimesHoHo_\alpha G),
\]
which is the usual Takesaki-Green induction process for crossed product $C^*$-algebras. The ideas of Proposition~\ref{genimpthm} show that a representation $\pi\rtimesHoHo U$ of $A\rtimesHoHo_\alpha G$ is induced from a representation $\tau\rtimes V$ of $A\rtimesHoHo_\alpha H$ if and only if there is a representation $\mu$ of $C_0(G/H)$ on $\H_\pi$ such that $(\pi\otimes \mu,U)$ is a covariant representation of $(A\otimes C_0(G/H),G,\alpha\otimes \lt)$; or, equivalently, such that every $\mu(f)$ commutes with every $\pi(a)$ and $(\mu,U)$ is a covariant representation of $(C_0(G/H),G,\lt)$.

One common problem with this general approach to imprimitivity theorems is to find a workable formula for the left action of an appropriate dense subspace on the right-Hilbert bimodule. In the case of Green's bimodule $X_H^G$, looking at the first equation in \cite[(B.5)]{BE} suggests that $z\in C_c(G,A)\subset A\rtimesHoHo_\alpha G$ should act on $x\in C_c(G,A)\subset X_H^G$ according to the formula
\begin{equation}\label{defleftactionX}
(z\cdot x)(s)=\int_G z(t)\alpha_t(x(t^{-1}s))\Delta_G(t)^{1/2}\,dt.
\end{equation}
However, $z$ is really a multiplier  of $(A\otimes C_0(G/H))\rtimesHoHo_{\alpha\otimes\lt}G$, and hence $z\cdot x:=(i_A\rtimesHeHe i_G(z))x$ is defined on elements of the form $x=b\cdot y$ by $z\cdot(b\cdot y)=(zb)\cdot y$. To verify that the left action of $z\in C_c(G,A)$ is indeed given by \eqref{defleftactionX}, we proceed as follows. First, we verify that for $b\in C_c(G\times G/H,A)$ and $y\in C_c(G,A)$ we have
\begin{equation}\label{Okonbx}
z\cdot(b\cdot y)=\int_G z(t)\alpha_t(b\cdot y(t^{-1}s))\Delta_G(t)^{1/2}\,dt,
\end{equation}
so that \eqref{defleftactionX} works for $x$ of the form $b\cdot y$.
Next we define a pairing $(z,x)\mapsto z*x$ using the formula on the right-hand side of \eqref{defleftactionX}, and check that
\begin{equation}\label{staradjoint}
\langle z*x,y\rangle_{A\rtimes H}=\langle x,z^* *y\rangle_{A\rtimes H}\ \mbox{ for $x,y\in C_c(G,A)\subset X$},
\end{equation}
where $z^*$ denotes the adjoint of $z$ in $C_c(G,A)\subset A\rtimesHoHo_\alpha G$. Then 
\begin{align*}
\langle (i_A\rtimesHeHe i_G(z))x,b\cdot y\rangle_{A\rtimes H}
&=\langle x, (i_A\rtimesHeHe i_G(z^*))(b\cdot y)\rangle_{A\rtimes H}\\
&=\langle x,z^**(b\cdot y)\rangle_{A\rtimes H}\qquad \mbox{ (by \eqref{Okonbx})}\\
&=\langle z*x,b\cdot y\rangle_{A\rtimes H}\qquad \mbox{ (by \eqref{staradjoint}),}
\end{align*}
and the density of $\{b\cdot y:b\in C_c(G\times G/H,A),\;y\in C_c(G,A)\}$ in $X_H^G$ now implies that $(i_A\rtimesHeHe i_G(z))x=z*x$ for every $x\in C_c(G,A)$. Thus the left action is given on the dense subspaces of compactly supported functions by \eqref{defleftactionX}, as claimed.

\section{The extension problem for actions}\label{sec-prob1}

Let $\alpha$ be an action of a locally compact group $G$ on a $C^*$-algebra $A$.
When $N$ is a closed normal subgroup of $G$, the action $\id\otimes\rt$ of $G/N$ on $A\otimes C_0(G/N)$ commmutes with the action $\alpha\otimes\lt$, and hence induces an action of $G/N$ on $(A\otimes C_0(G/N))\rtimesHoHo_{\alpha\otimes\lt} G$; we denote this induced action by $\beta$. Our solution to Problem~\ref{prob1} involves this action $\beta$.

\begin{thm}\label{thm-prob1}
Suppose that $\alpha$ is an action of a locally compact group $G$ on a $C^*$-algebra $A$, and that $N$ and $K$ are closed subgroups of  $G$ such that $N$ is normal in $G$ and $N\subset K$. Let $(\pi,U)$ be a covariant representation of $(A,N,\alpha)$ on $\H$.   Then there exists a covariant representation $(\pi,V)$ of $(A,K,\alpha)$ on $\H$ with $V|_N=U$ if and only if there is a representation $T$ of $K/N$ on $X_N^G\otimes_{A\rtimesHoHo_{\alpha} N}\H$ such that $(X_N^G\dashind(\pi\rtimesHoHo U), T)$ is a covariant representation of $((A\otimes C_0(G/N))\rtimesHoHo_{\alpha\otimes\lt}G, K/N,\beta)$.
\end{thm}

\begin{remark}
If we have realised $X_N^G\dashind(\pi\rtimesHoHo U)$ as the integrated form $(\rho\otimes \mu)\rtimes W$, then we can alternatively describe the covariance condition on $(X_N^G\dashind(\pi\rtimesHoHo U), T)$ as saying that every $T_{kN}$ commutes with every $\rho(a)$ and every $W_t$, and that $(\mu,T)$ is a covariant representation of $(C_0(G/N),K/N,\rt)$.
\end{remark}

As foreshadowed in the introduction, the main ingredients in the proof of Theorem~\ref{thm-prob1} will be a commutative diagram relating induction and restriction of representations and an imprimitivity theorem.  Recall that $(A\otimes C_0(G/N))\rtimesHoHo_{\alpha\otimes\lt}G$ is naturally isomorphic to $(A\rtimesHoHo_\alpha G)\rtimesHoHo_{\widehat\alpha|}(G/N)$ by \cite[Theorem~A.64]{BE}. When $G$ is amenable, so that $A\rtimesHoHo_{\alpha} G=A\rtimesHoHo_{\alpha,r} G$, and when $K=M$ is normal, so that it makes sense to restrict the coaction to $G/M$, \cite[Theorem~5.16]{BE} says that the following diagram of right-Hilbert bimodules commutes:
\begin{equation*}
\begin{CD}
(A\rtimesHoHo_\alpha G)\rtimesHoHo_{\widehat\alpha|}(G/N) @>\,\,X_N^G\,\,>>A\rtimesHoHo_{\alpha}N\\
@V{Y_{G/M}^{G/N}}VV @VV{A\rtimesHoHo_\alpha M}V\\
(A\rtimesHoHo_\alpha G)\rtimesHoHo_{\widehat\alpha|}(G/M)@>\,\,X_M^G\,\,>> A\rtimesHoHo_\alpha M.
\end{CD}
\end{equation*} 
From this, we deduce the commutativity of the following diagram relating induction and restriction of representations:
\begin{equation*}
\begin{CD}
\Rep((A\rtimesHoHo_\alpha G)\rtimesHoHo_{\widehat\alpha|}(G/N)) @<\,\,X_N^G\dashind\,\,<\cong<\Rep(A\rtimesHoHo_{\alpha}N)\\
@A{Y\dashind_{G/M}^{G/N}}AA @AA{\Res}A\\
 \Rep((A\rtimesHoHo_\alpha G)\rtimesHoHo_{\widehat\alpha|}(G/M))@<\,\,X_M^G\dashind\,\,<\cong< \Rep(A\rtimesHoHo_\alpha M).
\end{CD}
\end{equation*}
When $G$ is not amenable, the diagram in \cite[Theorem~5.16]{BE} involves the reduced crossed product $A\rtimesHoHo_{\alpha,r} G$, and hence is not useful here. So we need an analogue of \cite[Theorem~5.16]{BE} which uses the full crossed product $A\rtimesHoHo_\alpha G$ and allows $M$ to be non-normal. 
So suppose that $K$ is a closed subgroup of $G$ with $N\subset K$. Our analogue of \cite[Theorem~5.16]{BE} uses in place of Mansfield's bimodule $Y_{G/M}^{G/N}$ an imprimitivity bimodule $Z_{G/K}^{G/N}$ which we  now construct using the symmetric imprimitivity theorem of \cite{rae-sit}. 
\begin{prop} \label{prop-PPL}
Let $\alpha$ be an action of a locally compact group $G$ on a $C^*$-algebra $A$, let $K$ be a closed subgroup of $G$, and let $N$ be a closed normal subgroup of $G$ with $N\subset K$. To simplify the notation, set
\begin{equation}\label{defLR}
L:=((A\otimes C_0(G/N))\rtimesHoHo_{\alpha\otimes\lt}G)\rtimes_\beta K/N\mbox{ and\ }
R:=(A\otimes C_0(G/K))\rtimesHoHo_{\alpha\otimes\lt} G.
\end{equation} Then
$C_c(G\times G/N, A)$ completes to give an $L$--$R$ imprimitivity bimodule $Z_{G/K}^{G/N}$.
\end{prop}


\begin{proof}
This proof is similar to the proof of \cite[Proposition~1.1]{EKR}, so we omit some details.
Define left and right actions of $K/N\times G$ and
$G$ on $G/N\times G$ by \begin{equation*}
(kN,t)\cdot(sN,r)=(ksN, tr)\quad\text{and}\quad (sN,r)\cdot t=(stN, rt)
\end{equation*}
for $k\in K$ and  $s, t,r\in G$.  These actions are free and proper and
commute. Define   $\sigma:K/N\times G\to\Aut A$  by $\sigma_{(kN,t)}=\alpha_t$, and take
$\tau:G\to\Aut A$ to be the trivial action; clearly $\sigma$ and $\tau$ commute.  There are actions
$\kappa=\sigma\otimes\lt:K/N\times G\to\Aut(\Ind\tau)$ and $\omega=\tau\otimes\rt:G\to\Aut(\Ind\sigma)$ on the induced algebras $\Ind\tau$ and $\Ind\sigma$ such
that, for $f\in\Ind\tau$ and $g\in\Ind\sigma$,
\begin{gather*}
\kappa_{(kN,t)}(f)(sN,r)=\sigma_{(kN,t)}\big(f( (kN,t)^{-1}\cdot (sN,r)
)\big)=\alpha_t(f(k^{-1}sN, t^{-1}r)),\text{ and}\\
\omega_t(g)(sN, r)=\tau_t\big(g((sN,r)\cdot t)\big)=g(stN,rt).
\end{gather*}
By the symmetric imprimitivity theorem \cite[Theorem~1.1]{rae-sit},
$C_c(G/N\times G, A)$ completes to an $((\Ind\tau)\rtimes_\kappa (K/N\times
G))$-$((\Ind\sigma)\rtimes_\omega G)$ imprimitivity bimodule.

The isomorphisms $\Theta:C_0(G/N, A)\to\Ind\tau$  and $\Omega:C_0(G/K, A)\to\Ind\sigma$ given by
\begin{align*}
&\Theta(f)(sN,r)=f(rs^{-1}N),\quad \Theta^{-1}(g)(sN)=g(s^{-1}N, e),\\
&\Omega(f)(sN,r)=\alpha_r(f(s^{-1}K)),\quad \Omega^{-1}(g)(sK)=g(s^{-1}N, e)
\end{align*}
are suitably equivariant, hence induce isomorphisms
\begin{align*}
L:=(C_0(G/N, A)\rtimesHoHo_{\alpha\otimes\lt}G)\rtimes_\beta K/N&\cong (\Ind\tau)\rtimes_\kappa (K/N\times
G)\mbox{ and}\\
R:=C_0(G/K, A)\rtimesHoHo_{\alpha\otimes\lt} G&\cong(\Ind\sigma)\rtimes_\omega G.
\end{align*} At this stage we have an $L$--$R$ imprimitivity bimodule based on $C_c(G/N\times G, A)$,
and chasing through the construction and isomorphisms gives the  following formulas for the
actions and inner products:
\begin{gather*}
c\cdot y(sN, r)= \int_{K/N} \int_G c(kN, t,
rs^{-1}N)\alpha_t(y(k^{-1}sN,
t^{-1}r))\Delta_G(t)^{1/2}\Delta_{K/N}(kN)^{1/2}\, dt\, d(kN)\notag
\\
y\cdot b(sN, r)=\int_G y(st^{-1}N, rt^{-1})\alpha_{rt^{-1}}(b(t,
ts^{-1}))\Delta_G(t)^{-1/2}\, dt  \notag
\\
{}_{L}\langle x\,,\, y\rangle(kN, r, sN)
=\Delta_{K/N}(kN)^{-1/2}\Delta_G(r)^{-1/2}\int_G x(tN, st)\alpha_r(y(k^{-1}tN,
r^{-1}st)^*)\, dt  \notag
\\
\langle x\,,\, y\rangle_{R}(r, sK)
=\Delta_G(r)^{-1/2}\int_{K/N}\int_G \alpha_t\big( x(k^{-1}s^{-1}N,
t^{-1})^*y(k^{-1}s^{-1}rN, t^{-1}r) \big)\, dt\, d(kN)
 \end{gather*}
for $c\in C_c(K/N\times
G\times G/N, A)\subset L$, and
$b\in C_c(G\times G/K, A)\subset R$ and $x,y\in C_c(G/N\times G, A)$.
Following \cite{EKR}\footnote{The formula for the automorphism
$\Upsilon$ on the last line of \cite[page~156]{EKR} should be
$\Upsilon(x)(r,s)=x(s,sr^{-1})\Delta_G(s)^{1/2}$.}, we combine the inner products and actions with the vector space isomorphism $\Upsilon:C_c(G/N\times G)\to C_c(G\times G/N)$ defined by $\Upsilon(x)(sN,r)=x(r, rs^{-1}N)\Delta_G(r)^{1/2}$ to obtain our bimodule $Z_{G/K}^{G/N}$.
(We apply $\Upsilon$ to make our construction compatible with others, so that, for example, when $K=N=\{e\}$ and $B:=(A\otimes C_0(G))\rtimesHoHo_{\alpha\otimes\lt}G$ we recover the standard bimodule ${}_BB_B$; it is then easier to compare $Z_{G/K}^{G/N}$ with bimodules implementing Manfield imprimitivity.)
We end up with an $L$--$R$ imprimitivity bimodule based on $C_c(G\times G/N, A)$ with  the formulas 
\begin{align*}
c\cdot x(r,sN)&=\int_{K/N}\int_G c(kN,t,sN)\alpha_t(x(t^{-1}r,t^{-1}skN))
\Delta_{K/N}(kN)^{1/2}\, dt\, d(kN)\\
x\cdot b(r,sN)&=\int_G x(t,sN)\alpha_t(b(t^{-1}r,t^{-1}sK))\, dt\\
{}_L\langle x\,,\, y\rangle(kN,r,sN)&=
\Delta_{K/N}(kN)^{-1/2}\int_G x(t,sN)\alpha_r(y(r^{-1}t,r^{-1}skN)^*)\Delta_G(r^{-1}t)\, dt\\
\langle x\,,\, y\rangle_R(r, sK)&=\int_{K/N}
\int_G\alpha_t\big(x(t^{-1},t^{-1}skN)^*y(t^{-1}r,t^{-1}skN)  \big)\Delta_G(t)^{-1}
\, dt\, d(kN)
\end{align*}
for $c\in C_c(K/N\times G\times G/N, A)\subset L$, $b\in C_c(G\times G/K, A)\subset R$ and $x,y\in C_c(G\times G/N, A)$.
Note that when $N=\{e\}$  these formulas reduce to those of \cite[Proposition~1.1]{EKR}.
\end{proof}

Applying the construction of \S\ref{sec-genimp} to $Z_{G/K}^{G/N}$ and the canonical embedding $i_C$ of 
\begin{equation}\label{defC}
C:=(A\otimes C_0(G/N))\rtimesHoHo_{\alpha\otimes\lt}G
\end{equation}
in $M(L)=M(C\rtimes_\beta(K/N))$ gives a right-Hilbert $C$--$R$ bimodule and a corresponding induction process $Z\dashind_{G/K}^{G/N}$. Proposition~\ref{genimpthm} provides us with a ready-made imprimitivity theorem for this induction process:

\begin{prop}\label{impforZ}
Suppose $\pi$ is a representation of $C$ on $\H$. Then there is a representation $\tau$ of $R$ such that $\pi$ is equivalent to $Z\dashind_{G/K}^{G/N}\tau$ if and only if there is a unitary representation $V$ of $K/N$ on $\H$ such that $(\pi,V)$ is a covariant representation of $(C,K/N,\beta)$.
\end{prop}

The following commutative diagram is the analogue for $Z_{G/K}^{G/N}$ of \cite[Theorem~5.16]{BE}.

\begin{prop}\label{prop-cd} 
Suppose that $\alpha$ is an action of a locally compact group $G$ on a $C^*$-algebra $A$, and $N$ and $K$ are closed subgroups of  $G$ such that $N$ is normal in $G$ and $N\subset K$. Then the diagram
\begin{equation}\label{indresforZ}
\begin{CD}
(A\otimes C_0(G/N))\rtimesHoHo_{\alpha\otimes\lt}G @>\,\,X_N^G\,\,>>A\rtimesHoHo_{\alpha}N\\
@V{Z_{G/K}^{G/N}}VV @VV{A\rtimesHoHo_\alpha K}V\\
 (A\otimes C_0(G/K))\rtimesHoHo_{\alpha\otimes\lt}G@>\,\,X_K^G\,\,>> A\rtimesHoHo_\alpha K
\end{CD}
\end{equation}
commutes.
\end{prop}

The inclusion $i_C$ of $C$ in $M(C\rtimes_\beta(K/N))=M(L)$ induces a left action of $C$ on $Z_{G/K}^{G/N}$. For the proof of Proposition~\ref{prop-cd}, it is important that we can describe the left action of the dense subalgebra $C_c(G\times G/N,A)$ of $C$ on the subspace $C_c(G\times G/N,A)$ of the module $Z_{G/K}^{G/N}$.

\begin{lemma}\label{leftactonZasexpected} Suppose $N$ and $K$ are closed subgroups of  $G$ such that $N$ is normal in $G$ and $N\subset K$.
Let $f\in C_c(G\times G/N,A)\subset C$ and let $w\in C_c(G\times G/N,A)\subset Z_{G/K}^{G/N}$. Then $i_C(f)\cdot w$ is the element of $Z_{G/K}^{G/N}$ given by the compactly supported function
\begin{equation}\label{leftactionZ}
(i_C(f)\cdot w)(r,sN)=\int_G f(t,sN)\alpha_t(w(t^{-1}r,t^{-1}sN))\,dt.
\end{equation} 
\end{lemma}

\begin{proof}
As in the discussion at the end of \S\ref{sec-genimp}, we first claim that this formula is correct on elements of the form $w=l\cdot z$ for $l\in C_c(K/N\times G\times G/N,A)$ and $z\in C_c(G\times G/N,A)$. Notice that the right-hand side of \eqref{leftactionZ} is the formula for the convolution product $(c,d)\mapsto c*d$ on the subalgebra $C_c(G\times G/N,A)$ of $C=(A\otimes C_0(G/N))\rtimesHoHo_{\alpha\otimes\lt}G$, which we know is associative. The left action of $f$ is defined by $i_C(f)\cdot(l\cdot z)=(i_C(f)l)\cdot z$, where $i_C(f)l$ is the product of the multiplier $i_C(f)$ with the element $l$ of $C\rtimes_\beta(K/N)$; viewing $l$ as an element of $C_c(K/N,C)$, this product is given by the function $kN\mapsto f*l(kN)$ in $C_c(K/N,C)$. Thus
\begin{align*}
i_C(f_\cdot(l\cdot z)&(r,sN)=\\
&\int_{K/N}\int_G (f*l(kN))(t,sN)\alpha_t(z(t^{-1}r,t^{-1}skN))\Delta_{K/N}(kN)^{1/2}\,dt\,d(kN).
\end{align*}
By writing $\rt_{kN}(z)$ for the function $(r,sN)\mapsto z(r,skN)$, and recognising the integral over $G$ as a convolution, we can deduce that
\[
i_C(f)\cdot(l\cdot z)(r,sN)=\int_{K/N}\big((f*l(kN))*\rt_{kN}(z)\big)(r,sN)\Delta_{K/N}(kN)^{1/2}\,d(kN).
\]
Now applying associativity of $*$ to the right-hand side gives
\[
\int_{K/N}\big(f*(l(kN)*\rt_{kN}(z))\big)(r,sN)\Delta_{K/N}(kN)^{1/2}\,d(kN);
\]
expanding out the two convolutions, applying Fubini's theorem, and reinserting the definition of $\rt_{kN}(z)$ converts this to
\[
\int f(t,sN)\alpha_t\Big(\!\int\!\!\!\int l(kN,u,t^{-1}sN)\alpha_u(z(u^{-1}t^{-1}r,u^{-1}t^{-1}skN))\Delta_{K/N}(kN)^{1/2}du\,d(kN)\!\Big)dt,
\]
which we can recognise as $f*(l\cdot z)(r,sN)$. Thus
\begin{equation}\label{dotis*}
i_C(f)\cdot(l\cdot z)=f*(l\cdot z),
\end{equation}
as claimed.

A messy but routine computation shows that
\begin{equation}\label{pull*across}
\langle f*w,z\rangle_R=\langle w, f^**z\rangle_R
\end{equation}
for $f,w,z\in C_c(G\times G/N,A)$, where $f^*$ denotes the adjoint of $f$ in $C$. Now \eqref{dotis*} and \eqref{pull*across} give
\begin{equation}\label{lastcomp}
\langle i_C(f)w,l\cdot z\rangle_R=\langle w,i_C(f^*)(l\cdot z)\rangle_R=\langle w,f^**(l\cdot z)\rangle_R=\langle f*w,l\cdot z\rangle_R.
\end{equation}
Since the elements of the form $l\cdot z$ are dense in $Z_{G/K}^{G/N}$ (this is a general property of imprimitivity bimodules), \eqref{lastcomp} implies that $i_C(f)\cdot w=f*w$, as required.
\end{proof}

The right action of $A\rtimesHoHo_\alpha K$ on the module $X_K^G$ extends to an action of the multiplier algebra $M(A\rtimesHoHo_\alpha K)$, and so the inclusion $i$ of $A\rtimesHoHo_\alpha N$ in $M(A\rtimes K)$ induces a right action of $A\rtimesHoHo_\alpha N$ on  $X_K^G$.  Again, it is important for the proof of Proposition~\ref{prop-cd} that we can describe the right action of the dense subalgebra $C_c(N,A)$ of $A\rtimesHoHo_\alpha N$ on the dense subspace $C_c(G,A)$ of $X_K^G$.

\begin{lemma}\label{lem-expectedrightaction}
Suppose that $N$ and $K$ are closed subgroups of  $G$ such that $N$ is a normal subgroup of $K$. Let $f\in C_c(N,A)\subset A\rtimesHoHo_\alpha N$, $x\in C_c(G,A)\subset X_K^G$.  Then $x\cdot i(f)$ is the element of $X_K^G$ given by the compactly supported function on $G$ defined by
\begin{equation}\label{eq-expectedrightaction}
(x\cdot i(f)) (t)=\int_N x(tn)\alpha_{tn}(f(n^{-1}))
\Delta_N(n)^{-1/2}\, dn.
\end{equation}
\end{lemma}
\begin{proof}
Let  $f\in C_c(N,A)$. For $g\in C_c(K,A)$,  $i(f)g$ is the function in $C_c(K,A)$ given by $(i(f)g)(k)=\int_N f(n)\alpha_n(g(n^{-1}k))\, dn$, and thus \[(gi(f))(k)=(i(f^*)g^*)^*(k)=\int_N g(kn)\alpha_{kn}(f(n^{-1}))\, dn.\]

Temporarily denote the right-hand side of Equation~\eqref{eq-expectedrightaction} by $x\diamond f(t)$ and continue to write $R$ for $(A\otimes C_0(G/K))\rtimesHoHo_{\alpha\otimes\lt}G$. It suffices to show ${}_R\langle x\diamond f\,,\, y\cdot g\rangle ={}_R\langle x\cdot i(f)\,,\, y\cdot g\rangle$ for all $y\in C_c(G,A)$ and $g\in C_c(A,K)$.  The first step is to note that
\[((y\cdot g)\cdot i(f)) (t)=(y\cdot (gi(f))) (t)=(y\cdot g)\diamond f(t);\]
 the second equality follows from routine caculations (it helps to recall at the end that $\Delta_K(n)=\Delta_N(n)$ for all $n\in N$ since $N$ is normal in $K$). Next,
\begin{align*}
{}_R\langle &x\cdot i(f)\,,\,y\cdot g\rangle (s, tK)={}_R\langle x\,,\, y\cdot(gi(f^*))\rangle(s, tK)\\
&=\Delta_G(s)^{-1/2}\int_K x(tk)\alpha_s\big( y\cdot(gi(f^*)) (s^{-1}tk)^*\big)\, dk\\
&=\Delta_G(s)^{-1/2}\int_K x(tk)\alpha_s\big( ((y\cdot g)\diamond f^*) (s^{-1}tk)\big)^*\, dk\\
\intertext{which, using the formula  for $(y\cdot g)\diamond f^*$ established above followed by more routine calculations}
&={}_R\langle x\diamond f\,,\, y\cdot g\rangle(s, tK)
\end{align*}
as required.
\end{proof}

\begin{proof}[Proof of Proposition~\ref{prop-cd}] 
We continue to use the notations $L$, $R$ and $C$ established in \eqref{defLR} and \eqref{defC}. We will prove that 
\begin{equation}\label{eq-rh-iso2}
X_N^G\otimes_{A\rtimes N}\widetilde{X_K^G}\cong Z_{G/K}^{G/N}
\end{equation}
as right-Hilbert $C$--$R$ bimodules. Given this, we have 
\begin{equation}\label{indressuff}
X_N^G\otimes_{A\rtimes N}\widetilde{X_K^G}\otimes_R X_K^G\cong Z_{G/K}^{G/N}\otimes_R X_K^G
\end{equation}
as right-Hilbert $C$--$(A\rtimesHoHo_\alpha K)$ bimodules; since $\widetilde{X_K^G}\otimes_R X_K^G\cong A\rtimesHoHo_\alpha K$, \eqref{indressuff} implies that
\begin{equation*}
X_N^G\otimes_{A\rtimes N}(A\rtimesHoHo_\alpha K)\cong Z_{G/K}^{G/N}\otimes_{(A\otimes C_0(G/K))\rtimes G}X_K^G 
\end{equation*} 
as right-Hilbert $C$--$(A\rtimesHoHo_\alpha K)$ bimodules, which says precisely that \eqref{indresforZ} commutes.

We define $\Psi:C_c(G,A)\otimes C_c(G,A)\to C_c(G\times G/N,A)$  using the left inner product on $X_N^G$, so that
\[\Psi(x\otimes\tilde y)(r,sN)={}_C\langle x\,,\, y\rangle(r, sN)=
\Delta_G(r)^{-1/2}\int_N x(sn)\alpha_r(y(r^{-1}sn)^*)\, dn.
\]
We shall show that $\Psi$ extends to an isomorphism implementing \eqref{eq-rh-iso2}.  

We first show that $\Psi$ is isometric.
Let $x,y,v,w\in C_c(G,A)$. Then
\begin{align*}
\langle\Psi&(x\otimes\tilde y)\,,\,\Psi(w\otimes\tilde v)\rangle_R(r,sK)\\
&=\int_{K/N}\int_G\alpha_t\big(\Psi(x\otimes\tilde y)(t^{-1},t^{-1}skN)^*\Psi(w\otimes\tilde v)(t^{-1}r,t^{-1}skN)\big)\Delta_G(t)^{-1}\, dt\, d(kN)\\
&=\int_{K/N}\int_G\alpha_t\big({}_C\langle x\,,\, y\rangle(t^{-1},t^{-1}skN)^*{}_C\langle w\,,\, v\rangle(t^{-1}r,t^{-1}skN)\big)\Delta_G(t)^{-1}\, dt\, d(kN)\\
&=\int_{K/N}\int_G\int_N\int_N y(skn)\alpha_t\big(x(t^{-1}skn)^*w(t^{-1}skm)\big)\\
&\qquad\qquad\qquad\qquad\qquad\qquad\qquad\alpha_r(v(r^{-1}skm)^*)\Delta_G(r)^{-1/2}\, dm\, dn\, dt\, d(kN)\\
\intertext{which, after an application of Fubini's theorem and combining one of the integrals over $N$ with the integral over $K/N$ to get an integral over $K$, is equal to}
&\int_K\int_G\int_N y(skn)\alpha_t\big(x(t^{-1}skn)^*w(t^{-1}sk)\big)\alpha_r(v(r^{-1}sk)^*)
\Delta_G(r)^{-1/2}\, dn\, dt\, dk.
\end{align*}
We next note that 
\begin{align*}\langle x\otimes\tilde y\,,\, w\otimes\tilde v\rangle_R
&=\big\langle\langle w\,,\, x\rangle_{A\rtimesHoHo_\alpha N}\cdot\tilde y\,,\, \tilde v\big\rangle_R\\
&=\big \langle({y\cdot\langle x\,,\, w\rangle}_{A\rtimesHoHo_\alpha N})\widetilde{\ },\tilde v\big\rangle_R\\
&={}_R\big\langle y\cdot \langle x\,,\, w\rangle_{A\rtimesHoHo_\alpha N},v\big\rangle
\end{align*}
so that
\begin{align*}
\langle & x\otimes\tilde y\,,\, w\otimes\tilde v\rangle_R(r,sK)
=\Delta_G(r)^{-1/2}\int_K (y\cdot \langle x\,,\, w\rangle_{A\rtimesHoHo_\alpha N})(sk)\alpha_r(v(r^{-1}sk)^*)\, dk.
\intertext{Now we expand the right action of $A\rtimesHoHo_\alpha N$ on $X_K^G$ using the formula of Lemma~\ref{lem-expectedrightaction} and evaluate the $A\rtimesHoHo_\alpha N$-valued inner product on $X_N^G$ to get}
\langle x\otimes\tilde y&\,,\, w\otimes\tilde v\rangle_R(r,sK)\\
&=\Delta_G(r)^{-1/2}\int_K\int_N y(skn)\alpha_{skn}\big(\langle x\,,\, w\rangle_{A\rtimesHoHo_\alpha N}(n^{-1})\big)\alpha_r(v(r^{-1}sk)^*)\Delta_N(n)^{-1/2}\, dn\, dk\\
&=\Delta_G(r)^{-1/2}\int_K\int_N\int_G y(skn)\alpha_{sknt}\big(x(t^{-1})^*w(t^{-1}n^{-1})\big)\alpha_r(v(r^{-1}sk)^*)\, dt\, dn\, dk,\\
\intertext{which, by the change of variable $sknt\mapsto t$, is equal to}
&=\int_K\int_G\int_N y(skn)\alpha_t\big(x(t^{-1}skn)^*w(t^{-1}sk)\big)\alpha_r(v(r^{-1}sk)^*)
\Delta_G(r)^{-1/2}\, dn\, dt\, dk.
\end{align*}
Comparing the formulas for the two inner products shows that $\Psi$ is isometric.  

To see that $\Psi$ has dense range in $Z_{G/K}^{G/N}=\overline{C_c(G\times G/N,A)}$, we fix $z\in C_c(G\times G/N,A)$ and $\epsilon >0$. Green proved in \cite[Lemma~1.2]{G} that there is an approximate identity for $C$ consisting of elements $c$ of the form
\begin{equation}\label{Gapproxid}
c=\sum_{i=1}^n{}_C\langle x_i,y_i\rangle,
\end{equation}
where $x_i, y_i\in C_c(G,A)$ and ${}_C\langle\cdot,\cdot\rangle$ denotes the $C$-valued inner product on $C_c(G,A)\subset X_N^G$ which we used to define $\Psi$. The image of this approximate identity under $i_C:C\to M(L)$ converges strictly to $1_{M(L)}$, and thus there exists $c$ of the form \eqref{Gapproxid} such that 
\begin{equation}\label{CZdense}
\|z-c\cdot z\|_R<\epsilon. 
\end{equation}
Now we recall from Lemma~\ref{leftactonZasexpected} that $c\cdot z$ is given by the convolution product $c*z$ in $C_c(G\times G/N,A)\subset C$, and compute in the  imprimitivity bimodule $X_N^G$:
\[
c\cdot z=c*z=\sum_{i=1}^n{}_C\langle x_i,y_i\rangle*z=\sum_{i=1}^n{}_C\langle x_i,z^*\cdot y_i\rangle.
\]
But this says precisely that
\[
c\cdot z=\sum_{i=1}^n\Psi(x_i\otimes(z^*\cdot y_i)^{\widetilde{\ }}\;),
\]
so \eqref{CZdense} says that the range of $\Psi$ is dense.

Since $\Psi$ is isometric for the $R$-valued inner product and has dense range, $\Psi$ is automatically equivariant for the right actions of $R$. So it remains to check that $\Psi$ is equivariant for the left actions of $C$. But if $f\in C_c(G\times G/N,A)\subset C$ and $x,y\in C_c(G\times G/N,A)$, then 
\[
\Psi(f\cdot(x\otimes\tilde y))=\Psi((f\cdot x)\otimes\tilde y)={}_C\langle f\cdot x\,,\, y\rangle=f*{}_C\langle  x\,,\, y\rangle,
\]
which by Lemma~\ref{leftactonZasexpected} is precisely $f\cdot\Psi(x\otimes\tilde y)$. Thus $\Psi$ extends to the desired isomorphism \eqref{eq-rh-iso2}, and we have proved Proposition~\ref{prop-cd}. 
\end{proof}

\begin{proof}[Proof of Theorem~\ref{thm-prob1}]
The right-Hilbert bimodules $Z_{G/K}^{G/N}$ and $A\rtimesHoHo_\alpha K$  appearing in the commutative diagram of Proposition~\ref{prop-cd} induce maps $Z\dashind_{G/K}^{G/N}$ and $\Res$ on representations such that
\begin{equation}\label{eq-cd}
\begin{CD}
\Rep((A\otimes C_0(G/N))\rtimesHoHo_{\alpha\otimes\lt}G) @<\,\,X_N^G\dashind\,\,<\cong<\Rep(A\rtimesHoHo_{\alpha}N)\\
@A{Z\dashind_{G/K}^{G/N}}AA @AA{\Res}A\\
 \Rep((A\otimes C_0(G/K))\rtimesHoHo_{\alpha\otimes\lt}G)@<\,\,X_K^G\dashind\,\,<\cong< \Rep(A\rtimesHoHo_\alpha K)
\end{CD}
\end{equation}
commutes.

We are asking whether $\pi\rtimesHoHo U$ is in the range of $\Res:\Rep(A\rtimesHoHo_\alpha K)\to \Rep(A\rtimesHoHo_\alpha N)$; since the horizontal arrows in the commutative diagram~\eqref{eq-cd} are bijections, this happens precisely when $X_N^G\dashind(\pi\rtimesHoHo U)$ is in the range of $Z\dashind_{G/K}^{G/N}$, which by Proposition~\ref{impforZ} is equivalent to the existence of the representation $T$.  
\end{proof}

\section{The dual problem}\label{sec-state2}

Let $\delta:B\to M(B\otimes C^*(G))$ be a coaction of a locally compact group $G$ on a $C^*$-algebra $B$. 
We denote by $u$ the universal representation of $G$ in $UM(C^*(G))$, and by $w_G$ the unitary element of $M(C_0(G)\otimes C^*(G))$ given by the strictly continuous function $s\mapsto u(s):G\to UM(C^*(G))$.
A covariant homomorphism of $(B,G,\delta)$ into a $C^*$-algebra $D$ is a pair $(\pi,\mu)$ of homomorphisms $\pi:B\to M(D)$ and $\mu:C_0(G)\to M(D)$ such that, for every $b\in B$,
\[
(\pi\otimes\id)\circ\delta(b)=(\mu\otimes \id)(w_G)(\pi(b)\otimes 1)(\pi\otimes\id)(w_G)^*\ \mbox{ in $M(D\otimes C^*(G))$;}
\]
when $M(D)=B(\H)$, we call $(\pi,\mu)$ a covariant representation.
We denote the crossed product of $(B,G,\delta)$ by $(B\rtimesHoHo_\delta G,j_B,j_G)$, as in \cite[\S A5]{BE}. Then
\[
(j_B,j_G):(B,C_0(G))\to M(B\rtimesHoHo_\delta G)
\]
is universal for covariant homomorphisms on $(B,G,\delta)$ \cite[Theorem~A.41]{BE}; we write $\pi\rtimes\mu$ for the homomorphism on $B\rtimesHoHo_\delta G$ corresponding to a covariant homomorphism $(\pi,\mu)$.   

If $N$ is a normal subgroup of $G$, we denote by $\delta|$ the restriction of $\delta$ to a coaction of $G/N$, as in \cite[Example~A.28]{BE}. Viewing functions on $G/N$ as functions on $G$ gives a non-degenerate homomorphism of $C_0(G/N)$ into $C_b(G)=M(C_0(G))$, and hence it makes sense to restrict non-degenerate homomorphisms $\mu:C_0(G)\to M(D)$ to $C_0(G/N)$: first form the strictly continuous extension $\bar\mu:M(C_0(G))\to M(D)$, and then define $\mu|$ to be the composition of $\bar\mu$ with the inclusion of $C_0(G/N)$ in $M(C_0(G))$. (As in \cite{BE}, our notation will ignore the step of extending to $M(C_0(G))$.) Restricting the canonical injection $j_G$ of $C_0(G)$ in $M(B\rtimesHoHo_\delta G)$ gives a covariant homomorphism $(j_B,j_G|)$ of $(B,G/N,\delta|)$; more generally, if $(\pi,\mu)$ is covariant for $(B,G,\delta)$, then $(\pi,\mu|)$ is covariant for $(B,G/N,\delta|)$, and  $\pi\rtimes\mu|=(\pi\rtimes\mu)\circ(j_B\rtimes j_G|)$.    
If $N$ and $M$ are normal subgroups with $N\subset M$, then we can equally well  restrict covariant homomorphisms of $(B,G/N,\delta|)$ to $(B,G/M,\delta|)$, and our dual problem concerns this restriction process.

\begin{problem}\label{prob2}
Suppose $\delta:B\to M(B\otimes C^*(G))$ is a coaction of a locally compact group, $M$ is a closed normal subgroup of $G$ and $(\pi, \mu)$ is a covariant representation of $(B, C_0(G/M), \delta|)$.  For which closed normal subgroups $N$ of $G$ with $N\subset M$ is there a covariant representation $(\pi, \nu)$ of $(B, G/N, \delta|)$ such that $\mu=\nu|$?
\end{problem}

Our solutions of this problem are slightly different when $\delta$ is a maximal coaction, in which case we recover the stabilisation $B\otimes\K$ from $B\rtimesHoHo_\delta G$ as the full crossed product $(B\rtimesHoHo_\delta G)\rtimesHoHo_{\widehat\delta} G$ by the dual coaction, and when $\delta$ is a Quigg-normal coaction, in which case the Katayama duality theorem says that $B\otimes\K$ is isomorphic to the reduced crossed product $(B\rtimesHoHo_\delta G)\rtimesHoHo_{\widehat\delta,r} G$.

The maximal coactions of \cite{EKQ, KQFM} are by definition coactions $\delta:B\to M(B\otimes C^*(G))$ such that the canonical surjection 
\[
\big(((\id\otimes\lambda^G)\circ \delta)\rtimes (1\otimes M)\big)\rtimes (1\otimes\rho^G):(B\rtimesHoHo_\delta G)\rtimesHoHo_{\widehat\delta}G\to B\otimes \K(L^2 (G))
\]
is an isomorphism.  The main examples of maximal coactions are constructed from dual coactions on full crossed products, and include the restrictions of dual coactions to quotients by normal subgroups (see \cite[\S7]{KQFM}).
When $\delta$ is maximal and $N$ is a closed normal subgroup of $G$, we denote by $Y_{G/N}^G(B)$ or $Y_{G/N}^G$ the $((B\rtimesHoHo_\delta G)\rtimesHoHo_{\widehat\delta}N)$--$(B\rtimesHoHo_{\delta|}(G/N))$ imprimitivity bimodule constructed in \cite[\S5]{KQFM} using duality. This bimodule is the one used in \S\ref{sectmax}. 

A coaction $\delta:B\to M(B\otimes C^*(G))$ is Quigg-normal if the canonical map $j_B$ of $B$ into $M(B\rtimesHoHo_\delta G)$ is injective. The main examples of Quigg-normal coactions are dual coactions on reduced crossed products \cite[Proposition~2.3]{Q}, and the restrictions of these dual actions to quotients by normal subgroups \cite[Lemma~3.2]{KQ}. When $\delta$ is Quigg-normal and $N$ is a closed normal subgroup of $G$, we denote by $Y_{G/N}^G(B)$ or $Y_{G/N}^G$ the $((B\rtimesHoHo_\delta G)\rtimesHoHo_{\widehat\delta,r}N)$--$(B\rtimesHoHo_{\delta|}(G/N))$ imprimitivity bimodule constructed in \cite{KQ} (we need to combine Lemma~3.2 and Corollary~3.4 of \cite{KQ}). We emphasise that this bimodule, which is the one used in \S\ref{sectnor}, is essentially that constructed by Mansfield in \cite{M}, adapted for full coactions and avoiding amenability hypotheses.

\section{The dual problem for maximal coactions}\label{sectmax}

\begin{thm}\label{thm-prob2maximal}
Suppose that $\delta:B\to M(B\otimes C^*(G))$ is a non-degenerate maximal coaction of a locally compact group $G$ on a $C^*$-algebra $B$, that $N$ and $M$ are closed normal subgroups of $G$ with $N\subset M$, and that $(\pi, \mu)$ is a covariant representation of $(B, G/M, \delta|)$ on a Hilbert space $\H$.  Denote by $(\eta,U)$ the covariant representation of $(B\rtimesHoHo_\delta G, M,\widehat\delta)$ such that $\eta\rtimes U=Y_{G/M}^G\dashind(\pi\rtimes\mu)$.  
Then there exists a covariant representation $(\pi, \nu)$ of $(B, G/N, \delta|)$ on $\H$ such that $\nu|_{C_0(G/M)}=\mu$ if and only if there exists a representation $\phi$ of $C_0(M/N)$ in the commutant of $\eta(B\rtimesHoHo_\delta G)$ such that $(\phi, U)$ is a covariant representation of $(C_0(M/N), M,\lt)$.
\end{thm}

The proof of Theorem~\ref{thm-prob2maximal} again has two main ingredients: the following commutative diagram, which is the analogue for maximal coactions of \cite[Corollary~5.14]{BE}, and  Green's imprimitivity theorem.

\begin{prop}\label{resindmax}
Suppose $\delta:B\to M(B\otimes C^*(G))$ is a non-degenerate maximal coaction of a locally compact group $G$ on a $C^*$-algebra $B$, and  $N$ and $M$ are closed normal subgroups of $G$ with $N\subset M$. Then the following diagram of right-Hilbert bimodules commutes:
\begin{equation}
\label{514}
\xymatrix{
B\rtimesHoHo_\delta G \rtimesHoHo_{\widehat\delta|}M
\ar[rr]^{Y_{G/M}^G(B)}
\ar[d]_{X_N^M(B\rtimesHoHo_\delta G)}
&&
B\rtimesHoHo_{\delta|}(G/M)
\ar[d]^{B\rtimesHoHo_{\delta|}(G/N)}
\\
B\rtimesHoHo_\delta G \rtimesHoHo_{\widehat\delta|}N
\ar[rr]^{Y_{G/N}^G(B)}
&&
B\rtimesHoHo_{\delta|}(G/N).
}
\end{equation}
\end{prop}

\begin{proof}
We are going to mimic the proof of \cite[Corollary~5.14]{BE}, and for this we need some notation. 
If $\alpha: G\to\Aut A$ is an action, then there is a canonical isomorphism of $(A\otimes C_0(G/M))\rtimesHoHo_{\alpha\otimes\lt}G$ onto
$(A\rtimesHoHo_\alpha G)\rtimesHoHo_{\widehat\alpha|}(G/M)$ \cite[Theorem~A.64]{BE}. Applying this isomorphism to the dual action $\widehat\delta$ allows us to view the Green bimodule $X_M^G(B\rtimesHoHo_\delta G)$ as a $(B\rtimesHoHo_\delta G\rtimesHoHo_{\widehat\delta}G\rtimesHoHo_{\epsilon|}(G/M))$--$(B\rtimesHoHo_\delta G\rtimesHoHo_{\widehat\delta} M)$ imprimitivity bimodule, where $\epsilon$ denotes the double dual coaction of $G$ on $(B\rtimesHoHo_\delta G)\rtimesHoHo_{\widehat\delta} G$.  We also recall from \cite[Proposition~4.2]{KQFM} that the Katayama bimodule $K=K(B)$ 
is a $(B\rtimesHoHo_\delta G \rtimesHoHo_{\widehat\delta} G)$--$B$ imprimitivity bimodule with an
$\epsilon$--$\delta$ compatible coaction $\delta_K$ of $G$ on $K$.

As in the proof of \cite[Corollary~5.14]{BE}, we consider the following diagram which has \eqref{514} as its outer square:
\begin{equation}
\label{514bigger}
\xymatrix{
B\rtimesHoHo_\delta G \rtimesHoHo_{\widehat\delta}M
\ar[rr]^{Y_{G/M}^G(B)}
\ar[dddd]_{X_N^M(B\rtimesHoHo_\delta G)}
&&
B\rtimesHoHo_{\delta|}(G/M)
\ar[dddd]^{B\rtimesHoHo_{\delta|}(G/N)}
\\
&
B\rtimesHoHo_\delta G \rtimesHoHo_{\widehat\delta} G
\rtimesHoHo_{\epsilon|}(G/M)
\ar[ul]|{X_M^G(B\rtimesHoHo_\delta G)}
\ar[ur]|{K(B)\rtimesHoHo_{\delta_K|}(G/M)}
\ar[dd]|{B\rtimesHoHo_\delta G \rtimesHoHo_{\widehat\delta} G
\rtimesHoHo_{\epsilon|}(G/N)}
&
\\
&&\\
&
B\rtimesHoHo_\delta G \rtimesHoHo_{\widehat\delta} G
\rtimesHoHo_{\epsilon|}(G/N)
\ar[dl]|{X_N^G(B\rtimesHoHo_\delta G)}
\ar[dr]|{K(B)\rtimesHoHo_{\delta_K|}(G/N)}
&
\\
B\rtimesHoHo_\delta G \rtimesHoHo_{\widehat\delta}N
\ar[rr]^{Y_{G/N}^G(B)}
&&
B\rtimesHoHo_{\delta|}(G/N).
}
\end{equation}

It is shown in \cite[Corollary~6.4]{KQFM}  that the upper and lower triangles in \eqref{514bigger} commute. Since the diagonal arrows are implemented by imprimitivity bimodules and therefore invertible, to prove commutativity
of the outside square it suffices to prove that  the 
left and right quadrilaterals commute.

The modules in the left-hand quadrilateral are Green bimodules, so it is again convenient to work with an arbitrary action $\alpha$ of $G$ on a $C^*$-algebra $A$, and later take $A=B\rtimesHoHo_\delta G$ and $\alpha=\widehat\delta$.

In proving \cite[Proposition~8]{G} (which is induction in stages for Takesaki-Green  induction), Green showed that the map defined in terms of the right action of $C_c(M,A)\subset A\rtimesHoHo_\alpha M$ on $C_c(G,A)\subset X_M^G$ by $x\otimes y\mapsto x\cdot y$ extends an isomorphism of right-Hilbert $(A\rtimesHoHo_\alpha G)$--$(A\rtimesHoHo_\alpha N)$ bimodules from $X_M^G\otimes_{A\rtimes M}X_N^M$ to $X_N^G$. The 
left action of $C_c(G,A)\subset A\rtimesHoHo_\alpha G$ on $C_c(G,A)\subset X_M^G$ is given by the formula for the convolution product in $A\rtimesHoHo_\alpha G$ (see \eqref{defleftactionX}), and the action of $C_0(G/M)\subset C_b(G)$ by pointwise multiplication combines with this left action to give the left action of 
\[
\K(X_M^G)=(A\otimes C_0(G/M))\rtimesHoHo_{\alpha\otimes\lt}G.
\]
The algebra $C_0(G/M)$ also acts on $X_N^G$ through the inclusion of $C_0(G/M)$ in $C_b(G/N)=M(C_0(G/N))$, and this action is also by pointwise multiplication. The map $x\otimes y\mapsto x\cdot y$ preserves these left actions of $C_0(G/M)$, and hence is an isomorphism of right-Hilbert $((A\otimes C_0(G/M))\rtimesHoHo_{\alpha\otimes\lt}G)$--$(A\rtimesHoHo_\alpha N)$ bimodules.

On the other hand, the imprimitivity bimodule isomorphism $c\otimes x\mapsto c\cdot x$ of 
\[
((A\otimes C_0(G/N))\rtimesHoHo_{\alpha\otimes\lt}G)\otimes_{(A\otimes C_0(G/N))\rtimes G}X_N^G
\]
onto $X_N^G$ also gives an isomorphism of right-Hilbert $((A\otimes C_0(G/M))\rtimesHoHo_{\alpha\otimes\lt}G)$--$(A\rtimesHoHo_\alpha N)$ bimodules. When we identify $(A\otimes C_0(G/M))\rtimesHoHo_{\alpha\otimes\lt}G$ with $(A\rtimesHoHo_\alpha G)\rtimesHoHo_{\widehat \alpha|}(G/M)$, as in \cite[Theorem~A.64]{BE}, these two isomorphisms give the upper and lower triangles in a commutative diagram
\begin{equation*}
\xymatrix{A\rtimesHoHo_\alpha M\ar[d]_-{X_M^N(A)}&(A\rtimesHoHo_\alpha G)\rtimesHoHo_{\widehat\alpha|}(G/M)
\ar[d]^-{(A\rtimesHoHo_\alpha G)\rtimesHoHo_{\widehat\alpha|}(G/N)}\ar@{.>}[dl]_-{X^G_N(A)}\ar[l]_-{X^G_M}\\
A\rtimesHoHo_\alpha N&\ar[l]^-{X_N^G(A)} (A\rtimesHoHo_\alpha G)\rtimesHoHo_{\widehat\alpha|}(G/N)}
\end{equation*}
of right-Hilbert bimodules. Taking $(A,\alpha)=(B\rtimesHoHo_\delta G,\widehat\delta)$ gives the left-hand quadrilateral of~\eqref{514bigger}.

To study the right-hand quadrilateral, we consider the coaction $\nu$ of $G$ on the linking algebra $L(K)$ of $K=K(B)$ such that there is a canonical isomorphism of $L(K)\rtimesHoHo_{\nu} G$ onto $L(K\rtimesHoHo_{\delta_K}G)$ (see \cite[Lemma~3.10]{BE}). Then the natural map of $L(K)\rtimesHoHo_{\nu|} (G/M)$ into $M(L(K)\rtimesHoHo_{\nu|} (G/N))$ gives a non-degenerate 
imprimitivity bimodule homomorphism
of $K\rtimesHoHo_{\delta_K|}(G/M)$ into
$M(K\rtimesHoHo_{\delta_K|}(G/N))$  whose coefficient maps are the natural maps of 
$B\rtimesHoHo_\delta G \rtimesHoHo_{\widehat\delta} G\rtimesHoHo_{\epsilon|}(G/M)$ into $M(B\rtimesHoHo_\delta G \rtimesHoHo_{\widehat\delta} G\rtimesHoHo_{\epsilon|}(G/N))$ 
and
$B\rtimesHoHo_{\delta|}(G/M)$ into  $M(B\rtimesHoHo_{\delta|}(G/N))$. Since these 
last two maps give the two vertical 
arrows in the right-hand quadrilateral of \eqref{514bigger},
it follows from \cite[Lemma~4.10]{BE}
that the quadrilateral commutes.
\end{proof}

\begin{proof}[Proof of Theorem~\ref{thm-prob2maximal}]
From Proposition~\ref{resindmax} we obtain a commutative diagram
\begin{equation*}
\begin{CD}
\Rep((B\rtimesHoHo_\delta G)\rtimesHoHo_{\widehat\delta} M) @<\,\,Y_{G/M}^G\dashind\,\,<\cong<\Rep(B\rtimesHoHo_{\delta|}(G/M))\\
@A{X\dashind_N^M}AA @AA{\Res}A\\
 \Rep((B\rtimesHoHo_\delta G)\rtimesHoHo_{\widehat\delta} N)@<\,\,Y_{G/N}^G\dashind\,\,<\cong< \Rep(B\rtimesHoHo_{\delta|}(G/N)).
\end{CD}
\end{equation*}
We deduce that there exists an appropriate representation $\nu$ if and only if the representation $\eta\rtimes U=Y_{G/M}^G\dashind(\pi\rtimes\mu)$ of $(B\rtimesHoHo_\delta G)\rtimesHoHo_{\widehat\delta} M$ is induced from a representation of $(B\rtimesHoHo_\delta G)\rtimesHoHo_{\widehat\delta} N$. The theorem now follows from Green's imprimitivity theorem \cite[Theorem~6]{G} (see \S\ref{sec-genimp}).
\end{proof}

\subsection*{The extension problem for non-normal subgroups}

We now discuss possible reformulations of Theorem~\ref{thm-prob2maximal} in which $M$ and $N$ are replaced by non-normal subgroups $K$ and $H$, respectively. This discussion will lead us to the conclusion that we should be considering Quigg-normal coactions, rather than maximal ones.

Let $\delta$ be a coaction of $G$ on $B$. When $H$ is a closed subgroup of $G$ which is not normal, we have to decide how to interpret the crossed product by the homogeneous space $G/H$.  When $G$ is discrete or $\delta$ is a dual coaction, there is some choice here (see \cite{EQ} and \cite[\S 2]{EKR}), but in general the only available candidate is the \emph{reduced crossed product} $B\rtimesHoHo_{\delta, r}(G/H)$ of $B$ by the homogeneous space $G/H$, which is by definition\begin{equation}\label{defredcp}
B\rtimesHoHo_{\delta, r}(G/H):=\clsp\left\{j_B(b)j_G|(f):b\in B, f\in C_0(G/H) \right\}\subset M(B\rtimesHoHo_\delta G).
\end{equation}
It is explained in \cite{EKR} why this closed span is a $C^*$-algebra, why we think of it as a crossed product by $G/H$, and why this crossed product is called reduced. One way in which this behaves like a reduced crossed product rather than a full one is the apparent absence of a universal property: even when $N$ is normal, so that we already have a crossed product $B\rtimesHoHo_{\delta|}(G/N)$, the map $j_B\rtimes j_G|$ of $B\rtimesHoHo_{\delta|}(G/N)$ onto $B\rtimesHoHo_{\delta, r}(G/N)$ need not be an isomorphism. (We discuss this absence of a universal property further at the end of~\S\ref{sectnor}.)

The absence of a universal property means that we cannot construct representations of $B\rtimesHoHo_{\delta, r}(G/H)$ from covariant pairs, and hence we cannot restrict covariant representations to $B\rtimesHoHo_{\delta, r}(G/H)$. So when the larger subgroup $K$ is not normal, we can only restrict representations to $B\rtimesHoHo_{\delta, r}(G/K)$ from a subalgebra $C$ of $M(B\rtimesHoHo_\delta G)$ with $B\rtimesHoHo_{\delta, r}(G/K)\subset M(C)$. Fortunately, this applies when $C$ is the reduced crossed product by another homogeneous space of $G$:

\begin{lemma}\label{resOKonrcp}
Suppose that $H$ and $K$ are closed subgroups of $G$ with $H\subset K$. Then the inclusion of $B\rtimesHoHo_{\delta, r}(G/H)$ in $M(B\rtimesHoHo_\delta G)$ is non-degenerate, and 
\[
B\rtimesHoHo_{\delta, r}(G/K)\subset M(B\rtimesHoHo_{\delta, r}(G/H)).
\]
\end{lemma}

\begin{proof}
The inclusion of $B\rtimesHoHo_{\delta, r}(G/H)$ in $M(B\rtimesHoHo_\delta G)$ is non-degenerate because both $j_B$ and $j_G|$ are non-degenerate homomorphisms. Hence $M(B\rtimesHoHo_{\delta, r}(G/H))$ is naturally embedded in $M(B\rtimesHoHo_\delta G)$. It is clear that elements $j_B(c)$ multiply $B\rtimesHoHo_{\delta, r}(G/H)$ on the left, and that if $g\in C_b(G/K)$ then $j_G|(g)$ multiplies $B\rtimesHoHo_{\delta, r}(G/H)$ on the right; since taking adjoints shows that the elements of the form $j_G|(f)j_B(b)$ also span $B\rtimesHoHo_{\delta, r}(G/H)$, $j_B(c)$ and $j_G|(g)$ multiply on the other side as well. 
\end{proof}

So when the larger subgroup $K$ is not normal, the only extension problem which makes sense involves two reduced crossed products. There is still the possibility of an interesting extension problem when the larger subgroup $M$ is normal and $H$ is not. To pose such a problem, we have to decide how to make sense of the restriction map 
\begin{equation}\label{eq-res-wd}
\Res:\Rep(B\rtimesHoHo_{\delta,r}(G/H))\to\Rep(B\rtimesHoHo_{\delta|}(G/M)).
\end{equation} 
If $\rho$ is a representation of $B\rtimesHoHo_{\delta, r}(G/H)$, then composing with the canonical maps $j_B$ and $j_G|_{C_0(G/H)}$ gives a pair of representations
\begin{equation}\label{eq-restriction-pair}
\pi=\overline{\rho}\circ j_B\quad\text{and}\quad\mu=\overline{\rho}\circ j_G|_{C_0(G/H)}
\end{equation}
of $B$ and $C_0(G/H)$, respectively; the restriction of $\rho$ should be the pair $(\pi,\mu|)$, where $\mu|$ is the usual restriction of $\overline{\mu}$ to $C_0(G/M)\subset C_b(G/H)=M(C_0(G/H))$.  However, it takes work to see that $(\pi,\mu|)$ is a covariant representation of $(B, G/M,\delta|)$.

\begin{lemma}\label{lem-steve-lost}
The homomorphism $j_B\rtimes j_G|_{C_0(G/M)}$ of $B\rtimesHoHo_{\delta|}(G/M)$ into $M(B\rtimesHoHo_\delta G)$ has range in $M(B\rtimesHoHo_{\delta,r}(G/H))$, and for every non-degenerate representation $\rho$ of \text{$B\rtimesHoHo_{\delta, r} (G/H)$}, the representation $\overline{\rho}\circ\left( j_B\rtimes j_G|_{C_0(G/M)} \right)$ is the integrated form of 
\begin{equation}\label{eq-steve-lost}
(\overline{\rho}\circ j_B, \overline{\rho}\circ( j_G|_{C_0(G/H)}) |_{C_0(G/M)}).
\end{equation}
In particular, \eqref{eq-steve-lost} is a covariant representation of $(B, G/M,\delta|)$.
\end{lemma}

\begin{proof}
Since the range of $j_B\rtimes j_G|_{C_0(G/M)}$ is $B\rtimesHoHo_{\delta,r} (G/M)$, Lemma~\ref{resOKonrcp} implies that the range lies in $M(B\rtimesHoHo_{\delta, r}(G/H))$. We denote by $j_B^{B\rtimes(G/M)}$ the canonical homomorphism of $B$ into $M(B\rtimesHoHo_{\delta|}(G/M))$ (to distinguish it from $j_B:B\to M(B\rtimesHoHo_{\delta, r}(G/H))$).

If $\rho$ is a non-degenerate representation of $B\rtimesHoHo_{\delta, r}(G/H)$, then 
$\tau:=\overline{\rho}\circ \left(j_B\rtimes j_G|_{C_0(G/M)}\right)$ is a non-degenerate representation of $B\rtimesHoHo_{\delta|}(G/M)$, and hence is the integrated form of 
\[
(\overline{\tau}\circ j_B^{B\rtimes (G/M)}, \tau\circ j_{G/M}).\]  Now
\begin{align*}
\tau\circ j_B^{B\rtimes (G/M)}&=\overline{\overline{\rho}\circ ( j_B\rtimes j_G|_{C_0(G/M)})}\circ j_B^{B\rtimes (G/M)}=\overline{\rho}\circ j_B,\mbox{ and}\\
\overline{\tau}\circ j_{G/M}&=\overline{\overline{\rho}\circ ( j_B\rtimes j_G|_{C_0(G/M)})}\circ j_{G/M}
=\overline{\rho}\circ j_G|_{C_0(G/M)},
\end{align*}
which since $j_G|_{C_0(G/M)}=(j_G|_{C_0(G/H)})|_{C_0(G/M)}$ is the required covariant pair.
\end{proof}
 
Lemma~\ref{lem-steve-lost} says, first, that we have a well-defined restriction map \eqref{eq-res-wd} which is implemented by the homomorphism
\[
\phi=j_B\rtimes j_G|_{C_0(G/M)}:B\rtimesHoHo_{\delta|}(G/M)\to M(B\rtimesHoHo_{\delta, r}(G/H))
\]
and takes a representation to the pair $(\pi,\mu)$ described in \eqref{eq-restriction-pair}. Since the range of $\phi$ is (by definition) $B\rtimesHoHo_{\delta, r}(G/M)$ and since $B\rtimesHoHo_{\delta, r}(G/M)$ is a subalgebra of $M(B\rtimesHoHo_{\delta,r}(G/H))$, Lemma~\ref{lem-steve-lost} also implies that this restriction map factorises as
\[
\Rep(B\rtimesHoHo_{\delta, r}(G/H))\to\Rep(B\rtimesHoHo_{\delta,r}(G/M))\to\Rep(B\rtimesHoHo_{\delta|}(G/M)).
\]
Identifying the range of the second map amounts to determining the difference between $B\rtimesHoHo_{\delta|}(G/M)$ and $B\rtimesHoHo_{\delta, r}(G/M)$,  which is an interesting problem in its own right. Modulo solving this problem, we are led once again to the extension problem for reduced crossed products by homogeneous spaces.

To finish off this discussion, we want to show that this extension problem for reduced crossed products by homomgeneous spaces is effectively a problem about Quigg-normal coactions. Indeed, Quigg proved that for every coaction $\delta$ on $B$, there is a Quigg-normal coaction $\delta^n$ on a quotient $B^n$ of $B$ such that $B\rtimesHoHo_\delta G$ is naturally isomorphic to $B^n\rtimesHoHo_{\delta^n}G$ (see \cite{Q} or \cite[\S A.7]{BE}). Since $B\rtimesHoHo_{\delta, r}(G/H)$ is by definition a subalgebra of $M(B\rtimesHoHo_\delta G)$, it is naturally isomorphic to $B^n\rtimesHoHo_{\delta^n, r}(G/H)$. Thus the extension problem for reduced crossed products which we have arrived at is equivalent to the extension problem for the Quigg-normal coaction $\delta^n$. 

So the next step is to study the extension problem for a Quigg-normal coaction. As we shall see, even when both subgroups are normal, there are substantial difficulties to be overcome. The non-normal case then raises several additional technical problems which are of independent interest, and which we plan to discuss in a future paper.

\section{The dual problem for Quigg-normal coactions}\label{sectnor}

In this section, $\delta:B\to M(B\otimes C^*(G))$ is a Quigg-normal coaction, and $Y_{G/N}^G=Y_{G/N}^G(B)$ denotes the Mansfield bimodule of \cite{KQ}. Since the crossed product in the left coefficient algebra $(B\rtimesHoHo_\delta G)\rtimesHoHo_{\widehat\delta,r} N$ is a reduced crossed product, we need to discuss reduced Green induction before we can proceed.

Let $\alpha$ be an action of a group $G$ on a $C^*$-algebra $A$, $H$ a closed subgroup of $G$, and $I$  the kernel of the map $A\rtimesHoHo_\alpha H\to A\rtimesHoHo_{\alpha, r} H$.  A theorem of Quigg and Spielberg \cite{QS} implies that the quotient imprimitivity bimodule $X_{H,r}^G:=X_H^G/(X_H^G\cdot I)$ implements a Morita equivalence between the reduced crossed products $(A\otimes C_0(G/H))\rtimesHoHo_{\alpha\otimes\lt, r}G$ and $A\rtimesHoHo_{\alpha, r}H$ (see \cite{aHR-QS} or \cite[\S B.1]{BE} for further details). On the other hand, it is shown in \cite[Lemma 2.5]{KQR}, for example, that $X\dashind_H^G I$ is the kernel of the quotient map of $A\rtimesHoHo_{\alpha}G$ onto $A\rtimesHoHo_{\alpha,r}G$, and hence we can view $X_{H,r}^G$ as a right-Hilbert $(A\rtimesHoHo_{\alpha,r}G)$--$(A\rtimesHoHo_{\alpha, r} H)$ bimodule.  
Thus there is a well-defined induction process $X\dashind_{H,r}^G$ from $\Rep(A\rtimesHoHo_{\alpha, r}H)$ to $\Rep(A\rtimesHoHo_{\alpha, r}G)$.

If $(\eta,V)$ is a covariant representation of $(A,G,\alpha)$ such that $\eta\rtimes V$ factors through a representation of the reduced crossed product, then we write $\eta\rtimes_r V$ for the corresponding representation of $A\rtimesHoHo_{\alpha, r}G$.

\begin{thm}\label{thm-prob2normal}
Suppose $\delta:B\to M(B\otimes C^*(G))$ is a non-degenerate Quigg-normal coaction of a locally compact group $G$ on a $C^*$-algebra $B$, that $N$ and $M$ are closed normal subgroups of $G$ with $N\subset M$, and that $(\pi, \mu)$ is a covariant representation of $(B, G/M, \delta|_{G/M})$ on a Hilbert space $\H$.  Denote by $(\eta,U)$ the covariant representation of $(B\rtimesHoHo_\delta G, M,\widehat\delta)$ such that $\eta\rtimes_r U=Y_{G/M}^G\dashind(\pi\rtimes\mu)$.  
Then there exists a covariant representation $(\pi, \nu)$ of $(B, G/N, \delta|_{G/N})$ on $\H$ such that $\nu|_{C_0(G/M)}=\mu$ if and only if there exists a representation $\phi$ of $C_0(M/N)$ in the commutant of $\eta(B\rtimesHoHo_\delta G)$ such that $(\phi, U)$ is a covariant representation of $(C_0(M/N), M,\lt)$.
\end{thm}

From \cite[Theorem~3.1]{KQR} (or from \cite[Corollary~5.14]{BE}, which is the same theorem with a different proof) we obtain a commutative diagram
\begin{equation}\label{resindfromBE}
\begin{CD}
\Rep((B\rtimesHoHo_\delta G)\rtimesHoHo_{\widehat\delta,r} M) @<\,\,Y_{G/M}^G\dashind\,\,<\cong<\Rep(B\rtimesHoHo_{\delta|}(G/M))\\
@A{X\dashind_{N,r}^M}AA @AA{\Res}A\\
 \Rep((B\rtimesHoHo_\delta G)\rtimesHoHo_{\widehat\delta,r} N)@<\,\,Y_{G/N}^G\dashind\,\,<\cong< \Rep(B\rtimesHoHo_{\delta|}(G/N)).
\end{CD}
\end{equation}
From \eqref{resindfromBE} we deduce that there exists an appropriate representation $\nu$ if and only if the representation $Y_{G/M}^G\dashind(\pi\rtimes\mu)$ of $(B\rtimesHoHo_\delta G)\rtimesHoHo_{\widehat\delta,r} M$ is induced from a representation of $(B\rtimesHoHo_\delta G)\rtimesHoHo_{\widehat\delta,r} N$. 

So to complete the proof of Theorem~\ref{thm-prob2normal} we need a version of Green's imprimitivity theorem which allows us to recognise representations of $A\rtimesHoHo_{\alpha, r}M$ which have been induced via $X\dashind_{N,r}^M$ from representations of $A\rtimesHoHo_{\alpha,r}N$.  Theorem~\ref{redimpthm} below is exactly what we need, and should be of independent interest. In particular, it is a nice and apparently non-trivial application of non-abelian duality.

\begin{thm}\label{redimpthm}
Suppose that $\alpha$ is an action of a locally compact group $G$ on a $C^*$-algebra $A$ and $N$ is a closed normal subgroup of $G$. Let  $\pi\rtimes_r U$ be a representation of the reduced crossed product $A\rtimesHoHo_{\alpha,r} G$. Then there is a representation $\rho\rtimes_r V$ of the reduced crossed product $A\rtimesHoHo_{\alpha,r} N$ such that $\pi\rtimes_r U$ is equivalent to $X\dashind_{N,r}^G(\rho\rtimes_r V)$ if and only if there is a non-degenerate representation $\phi$ of $C_0({G/N})$ in the commutant of $\pi(A)$ such that $(\phi,U)$ is a covariant representation of $(C_0(G/N),G,\lt)$.
\end{thm}

The Morita equivalence $X_N^G$ descends to give an equivalence between the reduced crossed products $(A\otimes C_0(G/N))\rtimesHoHo_{\alpha\otimes\lt, r}G$ and $A\rtimesHoHo_{\alpha, r}N$, so
if $\pi\rtimes_r U$ is equivalent to $X\dashind_{N,r}^G(\rho\rtimes_r V)$, then the left action of $C_0(G/N)\subset M((A\otimes C_0(G/N))\rtimesHoHo_{\alpha\otimes\lt,r}G)$ on $X_{N,r}^G$ gives the required representation of $C_0(G/N)$ on $\H(X_{N,r}^G\dashind(\rho\rtimes_r V))=X_{N,r}^G\otimes_{A\rtimes N}\H_\rho$. 

Now suppose there is a representation $\phi:C_0(G/N)\to\pi(A)'$ such that $(\phi,U)$ is covariant for $(C_0(G/N),G,\lt)$; we have to show the existence of a suitable representation $\rho\rtimes_r V$. Green's imprimitivity theorem \cite[Theorem~6]{G} implies that there is a representation $\rho\rtimes V$ of the full crossed product $A\rtimesHoHo_\alpha N$ such that $(\pi\otimes \phi)\rtimes U\sim X_N^G\dashind(\rho\rtimes V)$. It suffices to show that $\rho\rtimes V$ factors through the reduced crossed product, and thus the result follows from the following proposition.

\begin{prop}\label{repfactors}
Suppose $\alpha$ is an action of a locally compact group $G$ on a $C^*$-algebra $A$, $N$ is a closed normal subgroup of $G$, and $(\rho, V)$ is a covariant representation of $(A,N,\alpha)$. Let $(\pi\otimes \phi,U)$ be the covariant representation of $(A\otimes C_0(G/N),G,\alpha\otimes\lt)$ such that $X_N^G\dashind (\rho\rtimes V)=(\pi\otimes \phi)\rtimes U$. If $\pi\rtimesHoHo U$ factors through a representation $\pi\rtimes_r U$ of the reduced crossed product $A\rtimesHoHo_{\alpha,r}G$, then $\rho\rtimes V$ factors through a representation $\rho\rtimes_r V$ of the reduced crossed product $A\rtimesHoHo_{\alpha,r} N$.
\end{prop}

\begin{proof}
It follows from \cite[Proposition~A.63]{BE} that $(\pi\rtimesHoHo U,\phi)$ is a covariant representation of $(A\rtimes_\alpha G,G/N,\widehat\alpha|)$, in the sense that
\begin{equation}\label{fullcov}
\big((\pi\rtimesHoHo U)\otimes\id\big)\circ\widehat\alpha|(b)
=\Ad (\phi\otimes\id(w_{G/N}))\big((\pi\rtimesHoHo U(b))\otimes 1\big)
\end{equation}
in $M(\K(\H_\pi)\otimes C^*(G/N))$ for every $b\in A\rtimesHoHo_\alpha G$. 
Let $u:G\to M(C^*(G))$ be the universal representation of $G$ and $(i_A^r,i_G^r)$ the canonical maps of $(A,G)$ into $M(A\rtimesHoHo_{\alpha,r}G)$.
We now recall from \cite[Example~A.27]{BE} that if  $Q$ is the quotient map of $A\rtimesHoHo_\alpha G$ onto $A\rtimesHoHo_{\alpha,r} G$, then the integrated form of $(i_A^r\otimes 1,i_G^r\otimes u)$ factors through a coaction $\widehat\alpha^n$ of $G$ on $A\rtimesHoHo_{\alpha,r} G$ which is characterised by
\begin{equation}\label{charalphan}
Q\otimes\id(\widehat\alpha(b))=\widehat\alpha^n(Q(b))\ \mbox{ for $b\in A\rtimesHoHo_\alpha G$.}
\end{equation}
Hitting both sides of \eqref{charalphan} with ${\id}\otimes q$, where $q:C^*(G)\to C^*(G/N)$ is the integrated form of the quotient map of $G$ onto $G/N$, gives the middle equality in the identity
\begin{equation}\label{anotherequ}
Q\otimes\id(\widehat\alpha|(b))=Q\otimes q(\widehat\alpha(b))={\id}\otimes q(\widehat\alpha^n(Q(b)))
=\widehat\alpha^n|(Q(b))\ \mbox{ for $b\in A\rtimesHoHo_\alpha G$.}
\end{equation}
By hypothesis, we have $\pi\rtimesHoHo U=(\pi\rtimes_r U)\circ Q$, so we can put $b=Q(c)$ in \eqref{anotherequ}, and then plug \eqref{anotherequ} into  \eqref{fullcov} to obtain
\begin{equation}\label{normcov}
\big((\pi\rtimes_r U)\otimes\id\big)\circ\widehat\alpha^n|(c)
=\Ad (\phi\otimes\id(w_{G/N}))\big((\pi\rtimes_r U(c))\otimes 1\big)\ \mbox{ for $c\in A\rtimes_{\alpha,r} G$.}
\end{equation}
Thus $(\pi\rtimes_r U,\phi)$ is a covariant representation of $(A\rtimesHoHo_{\alpha,r} G,G/N,\widehat\alpha^n|)$, and there is a representation $(\pi\rtimes_r U)\rtimes\phi$ of $(A\rtimesHoHo_{\alpha}G)\rtimesHoHo_{\widehat\alpha^n|}(G/N)$ such that
\[
((\pi\rtimes_r U)\rtimes\phi)(j_{A\rtimes G}(c)j_G(f))=\pi\rtimes_r U(c)\phi(f)\ \mbox{ for $c\in A\rtimesHoHo_{\alpha,r} G$ and $f\in C_0(G/N)$.}
\]

Theorem~A.65 of \cite{BE} says that there is an isomorphism
\[
\Omega:(A\rtimesHoHo_{\alpha,r}G)\rtimesHoHo_{\widehat\alpha^n|}(G/N)\to (A\otimes C_0(G/N))\rtimesHoHo_{\alpha\otimes\lt, r}G
\]
such that
\begin{align*}
\Omega\circ j_{A\rtimes_r G}\circ i_A^r(a)&=i^r_{A\otimes C_0(G/N)}(a\otimes 1),\\
\Omega\circ j_{A\rtimes_r G}\circ i_G^r(s)&=i^r_{G}(s),\ \mbox{ and}\\
\Omega\circ j_{G}(f)&=i^r_{A\otimes C_0(G/N)}(1\otimes f).
\end{align*}
Thus the representation $\big((\pi\rtimes_r U)\rtimes\phi\big)\circ \Omega^{-1}$ of $(A\otimes C_0(G/N))\rtimesHoHo_{\alpha\otimes\lt, r}G$ satisfies
\[
\big((\pi\rtimes_r U)\rtimes\phi\big)\circ \Omega^{-1}\circ R=(\pi\otimes \phi)\rtimes U,
\]
where $R$ is the quotient map of $(A\otimes C_0(G/N))\rtimesHoHo_{\alpha\otimes\lt}G$ onto the reduced crossed product. So $X_N^G\dashind(\rho\rtimes V)=(\pi\otimes \phi)\rtimes U$ factors through the reduced crossed product. 

Since Green's bimodule $X_N^G$ descends to give an equivalence between the reduced crossed products $(A\otimes C_0(G/N))\rtimesHoHo_{\alpha\otimes\lt, r}G$ and $A\rtimesHoHo_{\alpha, r}N$ and $X_N^G\dashind(\rho\rtimes V)$ factors through the reduced crossed crossed product, so does $\rho\rtimes V$. 
\end{proof}

This completes the proof of Theorem~\ref{redimpthm}, and hence also that of Theorem~\ref{thm-prob2normal}.

It is tempting to conjecture that Theorem~\ref{redimpthm} holds with $N$ replaced by an arbitrary closed subgroup $H$. Certainly, it holds when $H$ is amenable: given $\phi$, the usual imprimitivity theorem implies that $\pi\rtimesHoHo U$ is induced from a representation of $A\rtimesHoHo_{\alpha}H=A\rtimesHoHo_{\alpha,r}H$. However, the following example shows that Proposition~\ref{repfactors} cannot be extended to arbitrary subgroups, so that our proof of Theorem~\ref{redimpthm} breaks down. While this does not resolve the question of whether Theorem~\ref{redimpthm} can be extended, the example does suggest that there will be limitations on the extent to which a satisfactory representation theory can be developed for crossed products by homogeneous spaces. Example~\ref{ex-QS} was previously used by Quigg and Spielberg to provide an example of a dynamical system involving a commutative algebra which is not hyponormal \cite{QS}.

\begin{example}\label{ex-QS}
Consider the semidirect product $G=\R^2\rtimes SL_2(\R)$, the subgroup $H=SL_2(\R)$, and the quasi-regular representation $\Ind_H^G 1$ on $L^2(G/H)$. The map $(x,A)\mapsto x$ induces a homeomorphism of $G/H$ onto $\R^2$, and then
\[
\big((\Ind_H^G 1)_{(x,A)}\xi\big)(y)=\xi(A^{-1}(y-x))\ \mbox{ for $\xi\in L^2(\R^2)$.}
\]
Because $H$ is not amenable, the representation $1$ does not factor through the reduced crossed product $C_r^*(H)$. However, we shall prove that $\Ind_H^G 1$ does factor through $C_r^*(G)$.

The usual Fourier transform on $L^2(\R^2)$ implements a unitary equivalence of  $\Ind_H^G 1$ with the representation $U$ of $G$ given by
\[
(U_{(x,A)}\xi)(y)=e^{-2\pi ix\cdot y}\xi(A^ty)\ \mbox{ for $\xi\in L^2(\R^2)$.}
\] 
The representation $M$ of $C_0(\R^2)$ by multiplication operators gives a covariant representation $(M,U)$ of the system $(C_0(\R^2),G,\alpha)$ in which $\alpha_{(x,A)}(f)(y)=f(A^ty)$. Since $\{(0,0)\}$ has measure zero, $L^2(\R^2)=L^2(\R^2\setminus\{(0,0)\})$, and $(M,U)$ is also a covariant representation of $(C_0(\R^2\setminus\{(0,0)\}),G,\alpha)$ on $L^2(\R^2\setminus\{(0,0)\})$. Let $e_1=(1,0)$, let 
\[
H_1=\{A\in H:(A^t)^{-1}e_1=e_1\}=\bigg\{
\begin{pmatrix}1&0\\c&1\end{pmatrix}:c\in\R\bigg\},
\]
let $G_1=\R^2\rtimes H_1$, and define $\psi:G/G_1\to  \R^2\setminus\{(0,0)\}$ by $\psi((x,A)G_1)=(A^t)^{-1}e_1$. Then $\psi^*:f\mapsto f\circ \psi$ is an isomorphism of $C_0(\R^2\setminus\{(0,0)\})$ onto $C_0(G/G_1)$ such that $\psi^*(\alpha_{(x,A)}(f))=\lt_{(x,A)}(\psi^*(f))$. Thus $(M\circ(\psi^*)^{-1},U)$ is a covariant representation of $(C_0(G/G_1),G,\lt)$, and hence by the usual imprimitivity theorem, $U$ is induced from a representation of $G_1$. Since $H_1$ is isomorphic to $\R$, $G_1=\R^2\rtimes H_1$ is amenable, and every representation of the form $\Ind_{G_1}^G V$ factors through $C_r^*(G)$; hence so does $U\sim \Ind_H^G 1$.
\end{example}

To see the relevance of this example to crossed products by homogeneous spaces, let $\delta$ be a Quigg-normal coaction of $G$ on $B$ and $H$ a closed subgroup of $G$, and consider the reduced crossed product $B\rtimesHoHo_{\delta,r} (G/H)$ defined in \eqref{defredcp}. When $H=N$ is normal, it follows from \cite[Lemma~3.2]{KQ} that $j_B\rtimes j_G|_{C_0(G/N)}$ is an isomorphism of the usual coaction crossed product $B\rtimesHoHo_{\delta|}(G/N)$ onto $B\rtimesHoHo_{\delta,r}(G/N)$, and therefore $B\rtimesHoHo_{\delta,r}(G/N)$ is universal for covariant representations of $(B,G/N,\delta|)$. 

One naturally asks whether $B\rtimesHoHo_{\delta,r} (G/H)$ has a universal property when $H$ is not normal, but Example~\ref{ex-QS} suggests that this is unlikely. To see why, we consider the coaction $\delta_G^n:C^*_r(G)\to M(C^*_r(G)\otimes C^*(G))$, which is the integrated form of the representation $\lambda^G\otimes u$ (see \cite[Proposition~A.61]{BE}). The observation of the previous paragraph translates into:

\begin{lemma}\label{univpropCGN}
Suppose $N$ is a closed normal subgroup of $G$. Then $C^*_r(G)\rtimesHoHo_{\delta_G^n,r}(G/N)$ is universal for pairs $(\pi,\mu)$ consisting of a representation $\pi$ of $C^*_r(G)$ and a representation $\mu$ of $C_0(G/N)$ such that $(\mu,\pi\circ\lambda^G)$ is a covariant representation of $(C_0(G/N),G,\lt)$.
\end{lemma}

\begin{proof}
Since $C^*_r(G)\rtimesHoHo_{\delta^n_G|}(G/N)$ is isomorphic to $C^*_r(G)\rtimesHoHo_{\delta_G^n,r}(G/N)$, the latter is universal for covariant representations of $(C^*_r(G),G/N,\delta_G^n|)$. But it follows from the argument in the proof of \cite[Proposition~A.57]{BE} that a pair $(\pi,\mu)$ is a covariant representation of $(C^*_r(G),G/N,\delta_G^n|)$ if and only if $(\mu,\pi\circ\lambda^G)$ is a covariant representation of $(C_0(G/N),G,\lt)$.
\end{proof}

Notice that the assertion in Lemma~\ref{univpropCGN} makes sense when $N$ is replaced by a non-normal subgroup $H$. However, if we take $G$ and $H$ as in Example~\ref{ex-QS}, then the assertion is false. To see this, note that the representation $\Ind_H^G 1$ of $C^*(G)$ has the form $\pi\circ\lambda^G$, and the pair $(M,\Ind_H^G1)$ is a covariant representation of $(C_0(G/H),G,\lt)$. Since $1$ does not factor through $C^*_r(H)$, the Quigg-Spielberg theorem (as in \cite[Corollary~3]{aHR-QS}, for example) says that $M\rtimes (\Ind_H^G1)$ cannot factor through a representation $\tau$ of the reduced crossed product $C_0(G/H)\rtimesHoHo_{\lt,r}G$. 

On the other hand, if there were a representation $\rho$ of $C^*_r(G)\rtimesHoHo_{\delta_G^n,r}(G/H)$ such that $\rho(j_{C^*_r(G)}(z)j_{G/H}(f))=\pi(z)M(f)$, and $\phi$ is the isomorphism of $C_0(G/H)\rtimesHoHo_{\lt,r}G$ onto $C^*_r(G)\rtimesHoHo_{\delta_G^n,r}(G/H)$ from \cite[Proposition~2.8]{EKR}, then $\tau=\rho\circ\phi$ would give such a representation $\tau$. So the pair $(\pi,M)$ satisfies the covariance condition of Lemma~\ref{univpropCGN}, but cannot give rise to a representation of 
$C_0(G/H)\rtimesHoHo_{\lt,r}G$.

\appendix

\section{Compatibility of Mansfield bimodules}

In each of Sections~\ref{sec-prob1}, \ref{sectmax} and \ref{sectnor} we have used different induction processes for crossed products by coactions. In this appendix, we show that the induction processes used in Sections~\ref{sec-prob1} and \ref{sectmax} are generalisations of the original induction process of Mansfield, as adapted to full coactions by Kaliszewski and Quigg in \cite{KQ} and used in \S\ref{sectnor}. For simplicity, we suppose that there is just one subgroup $N$ and that this subgroup is normal.

First we consider the bimodule $Z^{G}_{G/N}$ used in \S\ref{sec-prob1}. Our current assumption that there is only one subgroup involved means that $Z^{G}_{G/N}$ is the bimodule constructed in \cite[Proposition~1.1]{EKR}, and it is shown in \cite[Theorem~4.1]{EKR} that when $K$ is normal and amenable, there is an isomorphism of $Z_{G/N}^{G}$ onto the Mansfield bimodule $Y_{G/N}^{G}$ which is compatible with the isomorphisms
\begin{align}
((A\otimes C_0(G))\rtimesHoHo_{\alpha\otimes\lt}G)\rtimes_\beta N
&\cong ((A\rtimesHoHo_\alpha G)\rtimesHoHo_{\widehat \alpha}G)\rtimesHoHo_{\widehat{\widehat\alpha}}N\label{isoleft}\\
(A\otimes C_0(G/N))\rtimesHoHo_{\alpha\otimes\lt}G
&\cong (A\rtimesHoHo_\alpha G)\rtimesHoHo_{\widehat \alpha|}(G/N)\label{isoright}
\end{align}
of coefficient algebras. When $N$ is not amenable, $Y_{G/N}^{G}$ is a proper quotient of $Z_{G/N}^{G}$ which implements a Morita equivalence between reduced crossed products (see \cite[Theorem~4.1]{EKR} for details).

In \S\ref{sectmax}, we consider a maximal coaction $\delta$ of $G$ on $B$ and use the bimodule from \cite{KQFM}, which we denote in this appendix by $DY_{G/N}^G$. In \cite{KQFM}, Kaliszewski and Quigg begin by considering a dual coaction $(C,G,\epsilon)$, and use the isomorphisms \eqref{isoleft} and \eqref{isoright} of coefficent algebras to view $Z_{G/N}^G$ as a $((C\rtimes_\epsilon G)\rtimesHoHo_{\widehat\epsilon} N)$--$(C\rtimesHoHo_{\epsilon|}(G/N))$ imprimitivity bimodule. They then apply this construction with $\epsilon$ the double dual coaction of $G$ on $C:=(B\rtimesHoHo_\delta G)\rtimesHoHo_{\widehat\delta} G$, and $DY_{G/N}^G$  is by definition the
$((B\rtimesHoHo_\delta G)\rtimesHoHo_{\widehat\delta} N)$--$(B\rtimesHoHo_{\delta|}(G/N))$ imprimitivity bimodule  which makes the following diagram commute:
\begin{equation}\label{defKQM}
\xymatrix{
B\rtimesHoHo_\delta G \rtimesHoHo_{\widehat\delta}N
\ar[rr]^{DY_{G/N}^G}
\ar[d]_{K(B)\rtimesHoHo_{\delta_K} G \rtimesHoHo_{\widehat\delta_K}N}
&&
B\rtimesHoHo_{\delta|}(G/N)
\ar[d]^{K(B)\rtimesHoHo_{\delta_K|}(G/N)}
\\
C\rtimes_\epsilon G \rtimesHoHo_{\widehat\epsilon}N
\ar[rr]^{Z_{G/N}^G}
&&
C\rtimesHoHo_{\epsilon|}(G/N).
}
\end{equation}
(see the proof of \cite[Theorem~5.3]{KQFM}). They check that, when
$(B,G,\delta)$ is a dual coaction,  $DY_{G/N}^G$ is isomorphic to $Z_{G/N}^G$  \cite[Proposition~6.5]{KQFM}.

It is also possible to apply the construction of \cite{KQ} to a maximal coaction $(B,G,\delta)$, yielding a $((B\rtimesHoHo_\delta G)\rtimesHoHo_{\widehat\delta,r} N)$--$(B\rtimesHoHo_{\delta, r}(G/N))$ imprimitivity bimodule $MY_{G/N}^G$ which is a direct generalisation of Mansfield's bimodule.
When $\delta$ is also Quigg-normal, $j_B\rtimes j_G|$ is an isomorphism of $B\rtimesHoHo_{\delta}(G/N)$ onto $B\rtimesHoHo_{\delta,r}(G/N)$, and $MY_{G/N}^G$ is the Mansfield bimodule we used in \S\ref{sectnor}. In general, as we shall see in the next paragraph,  $MY_{G/N}^G$ is a quotient of $DY_{G/N}^G$ which implements an equivalence between reduced crossed products.

To see that $MY_{G/N}^G$ can be a proper quotient of $DY_{G/N}^G$, consider again a dual coaction $(C,G,\epsilon)$. In \cite[Theorem~4.1]{EKR}, it is shown that the Rieffel correspondence for $Z_{G/N}^G(C)$ matches up the kernel of the regular representation of $(C\rtimes_\epsilon G)\rtimesHoHo_{\widehat\epsilon} N$ with the kernel of $j_C\rtimes j_G|$, and that the resulting quotient imprimitivity bimodule is isomorphic to $MY_{G/N}^G(C)$. Since the vertical arrows in \eqref{defKQM} also match up the kernels of the regular representations, we deduce that there is a commutative diagram
\begin{equation}\label{commW}
\xymatrix{
B\rtimesHoHo_\delta G \rtimesHoHo_{\widehat\delta,r}N
\ar[rr]^{W}
\ar[d]_{K(B)\rtimesHoHo_{\delta_K} G \rtimesHoHo_{\widehat\delta_K,r}N}
&&
B\rtimesHoHo_{\delta,r}(G/N)
\ar[d]^{K(B)\rtimesHoHo_{\delta_K,r}(G/N)}
\\
C\rtimes_\epsilon G \rtimesHoHo_{\widehat\epsilon,r}N
\ar[rr]^{MY_{G/N}^G(C)}
&&
C\rtimesHoHo_{\epsilon,r}(G/N),
}
\end{equation} 
in which the bimodule $W$ on the top arrow is the quotient of $DY_{G/N}^G$ associated to the kernel of $j_B\rtimes j_G|$. On the other hand, we know from the naturality\footnote{Strictly speaking, this is only proved in \cite[Theorem~4.21]{BE} for normal coactions. However, here we only need to apply naturality to the morphism associated to the Katayama imprimitivity bimodule $K(B)$, and the first part of  the proof of \cite[Theorem~4.21]{BE} with $B\rtimesHoHo_{\delta|}(G/N)$ replaced by $B\rtimesHoHo_{\delta,r}(G/N)$  establishes that this is OK.} of $MY_{G/N}^G$ that there is a similar diagram in which the top arrow is the Mansfield bimodule $MY_{G/N}^G(B)$; from this and \eqref{commW}, we deduce that $MY_{G/N}^G(B)$ is isomorphic to the quotient bimodule~$W$.

 \end{document}